\documentclass[12pt]{elsarticle}
\pdfoutput=1

\biboptions{sort&compress}

\global\long\def\vct#1{\boldsymbol{#1}}
\global\long\def\tns#1{{\boldsymbol{#1}}}
\global\long\def\discret#1{{\boldsymbol{\mathsf{#1}}}}
\global\long\def\discrets#1{{\mathsf{#1}}}
\global\long\def\transpose#1{#1^{\mathsf{T}}}

\global\long\def\idx#1{{\mathsf{#1}}}

\def\solidsquare{\vrule height .9ex width .8ex depth -.1ex\nobreak\ }
\def\solidcopen{\drawline{10}{.5}\nobreak\raise 0.5pt\hbox{\circle}\drawline{10}{.5}\nobreak\ }

\usepackage{amssymb}
\usepackage{amsmath}
\usepackage{lineno}
\usepackage{booktabs}
\usepackage{bm}
\usepackage{subcaption}
\usepackage{xcolor}  

\usepackage{ctable}

\definecolor{color1}{HTML}{000000} 
\definecolor{color2}{HTML}{4169E1} 
\definecolor{color3}{HTML}{8B0000} 
\definecolor{color4}{HTML}{006400} 

\journal{Journal of Computational Physics}

\begin{document}

\begin{frontmatter}

\title{A cut-cell finite volume -- finite element coupling approach for fluid-structure 
interaction in compressible flow}


\author{Vito Pasquariello$^1$\corref{cor}}
\cortext[cor]{Corresponding author. Tel.: +49 89 289 16134; fax: +49 89 289 16139}
\ead{vito.pasquariello@tum.de}

\author{Georg Hammerl$^2$\corref{}}

\author{Felix \"Orley$^1$\corref{}}

\author{Stefan Hickel$^{1,3}$\corref{}}
       
\author{\\Caroline Danowski$^2$\corref{}}

\author{Alexander Popp$^2$\corref{}}

\author{Wolfgang A. Wall$^2$\corref{}}

\author{\\Nikolaus A. Adams$^1$\corref{}}

\address{$^1$Institute of Aerodynamics and Fluid Mechanics, Technische Universit\"at M\"unchen\\
         Boltzmannstr.~15, 85748 Garching bei M\"unchen, Germany}
         
\address{$^2$Institute for Computational Mechanics, Technische Universit\"at M\"unchen\\
         Boltzmannstr.~15, 85748 Garching bei M\"unchen, Germany}

\address{$^3$Aerodynamics Group, Faculty of Aerospace Engineering, Technische Universiteit Delft\\
         Kluyverweg 1, 2629 HS Delft, The Netherlands}

\begin{abstract}
  We present a loosely coupled approach for the solution of fluid-structure interaction problems 
between a compressible flow and a deformable structure. The method is based on staggered 
Dirichlet-Neumann partitioning. The interface motion in the Eulerian frame 
is accounted for by a conservative cut-cell Immersed Boundary method. The present approach enables sub-cell 
resolution by considering individual cut-elements within a single fluid cell, which guarantees an accurate 
representation of the time-varying solid interface. The cut-cell procedure inevitably leads to non-matching 
interfaces, 
demanding for a special treatment. A Mortar method is chosen in order to obtain a conservative and consistent 
load transfer.
We validate our method by investigating two-dimensional test cases comprising a shock-loaded rigid cylinder and a 
deformable panel. Moreover, the aeroelastic instability of a thin plate structure is studied with a focus on the 
prediction of flutter onset.
Finally, we propose a three-dimensional fluid-structure interaction test case of a flexible inflated 
thin shell interacting with a shock wave involving large and complex structural deformations.

\end{abstract}


\end{frontmatter} 


\section{Introduction}
\label{sec:intro}
Compressible fluid-structure interaction (FSI) occurs in a broad range of technical applications involving, e.g., 
nonlinear aeroelasticity \cite{Piperno1997,Farhat2000} and shock-induced deformations of rocket nozzles 
\cite{Garelli2010, Zhao2013}. The numerical modeling and simulation of compressible FSI can be challenging, in 
particular if an accurate representation of the structural interface within the fluid solver and a consistent coupling 
of both subdomains is required.

FSI algorithms are generally classified as monolithic or partitioned. 
One main advantage often attributed to monolithic approaches is their numerical robustness due to solving a single
system which includes the full information of the coupled nonlinear FSI problem. 
On the other hand, partitioned algorithms for FSI are often used because they facilitate the coupling
of different specialized single-field solvers. A further distinction can be made between loosely and 
strongly coupled algorithms, depending on whether the coupling conditions are satisfied exactly at each time 
step, or not. While partitioned algorithms can be made strong by introducing equilibrium iterations 
\cite{Kuettler2008}, loosely coupled approaches are more frequently used in the field of aeroelasticity and 
compressible flows in general \cite{Farhat2000, Cirak2007}. A disadvantage of loosely coupled partitioned algorithms is 
the artificial added mass effect \cite{Causin2005, Forster2007}, which may lead to numerical instability in 
incompressible flows and for high fluid-solid density ratios. Recently, 
so-called Added-Mass Partitioned algorithms have been developed for compressible fluids interacting with 
rigid and elastic solids \cite{Banks2013,Banks2012} as well as for incompressible fluids \cite{Banks2014a}. 
These methods allow to overcome the added mass instability by formulating appropriate fluid-structure interface 
conditions.

FSI problems involve a load and motion transfer at the conjoined interface. In the simple case of matching fluid 
and solid discretization, this results in a trivial task. However, different resolution requirements within the 
fluid and solid fields lead to non-matching discrete interfaces. 
An overview of existing coupling methods for 
non-matching meshes can be found in \cite{DeBoer2007}. Simple methods such as nearest-neighbor interpolation 
and projection methods are frequently used \cite{Farhat1998, joppich:mpcci}. 
The mentioned methods do not conserve angular momentum across the interface. 
Consistency can be achieved with more sophisticated approaches, such as weighted residual methods, which 
introduce Lagrange multipliers as additional interface variables. In this context, Mortar methods have first 
been proposed for non-overlapping domain decomposition in \cite{Bernardi1994}, enhanced with dual
shape functions for the Lagrange multipliers in \cite{Wohlmuth2000} and applied to FSI problems and 
mesh tying in fluid flow, e.g. in \cite{Kloeppel2011, Ehrl2014}. While Mortar methods introduce Lagrange 
multipliers only on one side of the interface, Localized Lagrange Multipliers consider them on both sides 
of the interface \cite{Ross2008}.

Another classification of FSI methods is based on the 
representation of the time-varying solid interface within the fluid domain. Two main approaches can be distinguished in 
this context, which are Arbitrary Lagrangian Eulerian (ALE) methods 
\cite{Donea:1982ab,Farhat1995}, and Immersed Boundary Methods (IBM) \cite{Peskin:1972ty, Mittal:2005ii}. ALE 
approaches employ body-fitted grids, hence requiring a mesh evolution algorithm. This task may 
be complex in case of large solid displacements. On the other hand, IBM often operate 
on fixed Cartesian fluid 
grids, making this type of approach very appealing for the simulation of flows past complex geometries and for the 
solution of FSI problems with large deformations.
IBM, such as continuous forcing and ghost-cell approaches, may suffer from spurious loss or 
production of mass, momentum and energy at the interface \cite{Mittal:2005ii}. Such non-conservativity poses a 
particular problem for large-eddy simulations, which employ coarse grids and rely on an accurate flow prediction in 
near-wall regions over large time scales. Moreover, the accurate capturing of shocks is 
based on conservation properties.
Conservativity is recovered with Cartesian cut-cell methods, which were first introduced by 
\citet{Clarke:1986wg} and \citet{GAFFNEYJR:2012en} for inviscid flows and later extended to viscous flows 
by \citet{Udaykumar:1996tj} and \citet{Ye:1999vu}. In this method, the finite volume 
cells at the boundaries are reshaped to fit locally the boundary surface with a sharp interface, which in turn assures 
strict conservation of mass, momentum and energy. A drawback of cut-cell methods is that the 
fluid volume fraction of cut-cells may become very small and therefore can lead to numerical instability with explicit 
time integration schemes. A stabilization of the underlying time integration scheme can be achieved by so-called 
cell-merging \cite{Ye:1999vu}, cell-linking \cite{Kirkpatrick:2003wr} or flux redistribution techniques 
\cite{Colella:2006kq,Hu2006}.

In this paper we develop a loosely coupled approach for the solution of FSI problems between a 
compressible fluid and a deformable structure. We employ the Finite Volume Method (FVM) for solving the Euler equations 
on Cartesian grids and the Finite Element Method (FEM) for solving the structural problem. The interface 
motion is accounted for by a conservative cut-cell IBM. Previous proposed methods reconstruct the interface 
geometry based on a level-set function \cite{Meyer:2010bg, grilli2009conservative, Gunther2014}. 
\citet{Oerley:2015} developed a conservative cut-element method that allows for representing the 
fluid-solid interface with sub-cell resolution for rigid body motion. We extend this method to arbitrary interface 
deformations.
The combination of a cut-element IBM with a Mortar method for coupling of the solid and fluid subdomains in a 
consistent and efficient way is the essential new contribution of this paper.

This paper is structured as follows: First, the governing equations for fluid and solid and the 
fluid-structure interface conditions are introduced in 
Section~\ref{sec:governing_equations}. Section~\ref{sec:numerical_approach_fluid} gives a detailed overview on the 
numerical treatment of moving 
boundaries together with the discretization methods used for the fluid. The FEM used to solve the 
structural problem is presented in 
Section~\ref{sec:str:fem}. In Section~\ref{sec:coupling_procedure}, the staggered coupling algorithm is 
presented together with the new coupling approach for non-matching interfaces. In 
Section~\ref{sec:validation}, the method is 
validated with well-established two-dimensional test cases and a convergence study is 
presented. In Section~\ref{sec:numerical_example}, 
we propose a new test case for the interaction between a flexible inflated thin shell and a shock wave, demonstrating 
in particular the capability of our FSI approach to handle large three-dimensional deformations. Concluding 
remarks are given in Section~\ref{sec:conclusions}.

\newpage
\section{Mathematical and physical model}
\label{sec:governing_equations}

\begin{figure}
  \centering
  \includegraphics[width= 0.6\columnwidth]{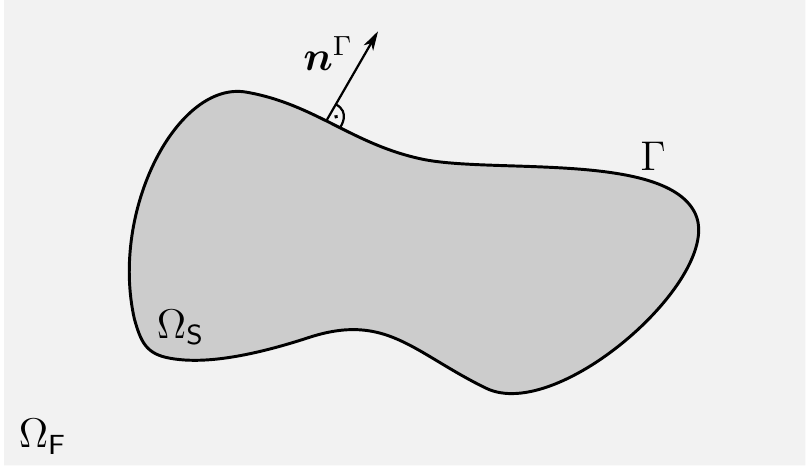}
  \caption{Schematic of FSI domain.}
  \label{fig:computational_domain}  
\end{figure}

As depicted in Fig.~\ref{fig:computational_domain}, the computational domain is divided into a fluid and solid domain, 
$\Omega_{\idx{F}}$ and 
$\Omega_{\idx{S}}$, respectively. The conjoined interface is denoted as $\Gamma =\Omega_{\idx{F}} \cap 
\Omega_{\idx{S}}$ and its normal vector 
$\vct{n}^{\Gamma}$ in spatial configuration points from the
solid into the fluid domain.

\subsection{Governing equations for the fluid}
\label{sec:governing_fluid}

We consider the three-dimensional, fully compressible Euler equations in conservative form
\begin{equation}
\frac{\partial\vct{w}}{\partial t} + \nabla\cdot\tns{K}(\vct{w})= 0\ {\rm in} \ \Omega_{\idx{F}} \,.
\label{eqn:nse}
\end{equation}
The state vector $\vct{w}=[\rho_\idx{F}, \rho_\idx{F} u_1, \rho_\idx{F} u_2, \rho_\idx{F} u_3,E_{\idx{t}}]$ 
contains the conserved variables density $\rho_\idx{F}$, momentum $\rho_\idx{F} \vct{u}$ and total energy 
$E_{\idx{t}}$. The subscript $\idx{F}$ denotes fluid quantities and is used whenever a distinction between both 
subdomains is necessary. The individual contributions of the flux tensor $\tns{K}=(\vct{f}, \vct{g}, 
\vct{h})$ are given as
\begin{equation}
\vct{f}(\vct{w})=
\begin{pmatrix}
\rho_\idx{F} u_1 \\
\rho_\idx{F} {u_1}^{2}+p \\
\rho_\idx{F} u_1u_2 \\
\rho_\idx{F} u_1u_3 \\
u_1(E_{\idx{t}}+p)
\end{pmatrix},\,
\vct{g}(\vct{w})=
\begin{pmatrix}
\rho_\idx{F} u_2\\
\rho_\idx{F} u_2u_1\\
\rho_\idx{F} {u_2}^{2}+p\\
\rho_\idx{F} u_2u_3\\
u_2(E_{\idx{t}}+p)
\end{pmatrix},\,
\vct{h}(\vct{w})=
\begin{pmatrix}
\rho_\idx{F} u_3\\
\rho_\idx{F} u_3u_1\\
\rho_\idx{F} u_3u_2\\
\rho_\idx{F} {u_3}^{2}+p\\
u_3(E_{\idx{t}}+p)
\end{pmatrix}\, ,
\label{eq:ConvectiveFluxes}
\end{equation}
where $p$ is the static pressure.
We consider a perfect gas with a specific heat ratio of $\gamma 
=1.4$ and specific gas constant of $\text{R} = 287.058\,\frac{\text{J}}{\text{kg}\cdot\text{K}}$.
The total energy is given by
\begin{equation}
E_{\idx{t}}=\frac{1}{\gamma-1}p+\frac{1}{2}\rho_\idx{F} u_i u_i\, ,
\label{eq:tot_energy}
\end{equation}
assuming an ideal gas equation of state $p = \rho_\idx{F} \text{R} T$, where $T$ is the static temperature. If not stated 
otherwise, we use the Einstein summation convention.

\subsection{Governing equations for the solid}
\label{sec:governing_solid}

The structural field is governed by the local form of the balance of linear momentum
\begin{equation}\label{equ:str:momRKloc} 
  \rho_{\idx{S};0} \, \ddot{\vct{d}}
  \,=\, \nabla_0 \cdot \,  ( \tns{F} \cdot \tns{S} ) \,+\, \hat{\vct{b}}_0 
  \ {\rm in} \ \Omega_{\idx{S}} \,,
\end{equation}
describing equilibrium of the forces of inertia, internal and external forces in the undeformed structural domain
$\Omega_{\idx{S}}$. Herein $\nabla_0 \cdot (\bullet)$ is the material divergence operator and the index $\idx{S}$
represents the domain of the structural problem. The structural material density is denoted by $\rho_{\idx{S};0}$.
Furthermore, $\vct{d}$
and $\ddot{\vct{d}}$ are the unknown displacements and accelerations, respectively. The vector field $\hat{\vct{b}}_0$
is the given material body force. The internal
forces are expressed in terms of the second Piola-Kirchhoff stress tensor $\tns{S}$ and the deformation gradient
$\tns{F}$. 

To determine the stresses, various constitutive laws can be used. For the sake 
of simplicity, in this work a
hyperelastic Saint Venant-Kirchhoff material model with strain energy density function $\Psi$ per unit reference volume
is chosen as
\begin{equation}\label{equ:strainenergy}
  \Psi(\tns{E})
  \,=\, \mu_\idx{S} \,\tns{E} : \tns{E}
  \,+\, \dfrac{1}{2} \, \lambda_\idx{S} \, ( \tns{E} : \tns{I} )^2
  \,,
\end{equation}
with the Lam\'{e} constants $\lambda_\idx{S}$ and $\mu_\idx{S}$ and the second-order identity tensor $\tns{I}$.
The Green-Lagrange strain tensor is defined as
\begin{equation}\label{equ:GLstrain}
  \tns{E} \,=\, \dfrac{1}{2} \, ( \tns{F}^{\sf T} \cdot \tns{F} \,-\, \tns{I} ) \,.
\end{equation}
The second Piola-Kirchhoff stress
\begin{equation}\label{equ:PK2stressconstr}
  \tns{S} \,=\, \dfrac{\partial \Psi}{\partial \tns{E}}
\end{equation}
is derived using \eqref{equ:strainenergy}.
Alternatively, the first Piola-Kirchhoff stress tensor
\begin{equation}
 \tns{P} = \tns{F} \cdot \tns{S}
\label{eq:cauchy_stress_solid}
\end{equation} 
may be used.

The boundary of the structural field $\partial \Omega_{\idx{S}}$ is divided into pairwise disjoint boundary segments
\begin{equation}
 \partial\Omega_{\idx{S}} \, = \, \Gamma_{\idx{S;D}} \cup \Gamma_{\idx{S;N}} \cup \Gamma \, .
\end{equation}
On the Dirichlet boundary $\Gamma_{\idx{S;D}}$, the displacements are prescribed, whereas on the Neumann boundary
$\Gamma_{\idx{S;N}}$, the traction vector $\hat{\vct{t}}_0$ is prescribed using the unit normal vector $\vct{n}_0$ 
in material configuration. Thus, the boundary conditions
\begin{align}
  \label{equ:str:dbc}
  \vct{d} \,&=\, \hat{\vct{d}} \quad {\rm on} \  \Gamma_{\idx{S;D}} \,,
  \\ 
  \label{equ:str:nbc}
  \tns{P} \cdot \vct{n}_0 \,&=\, \hat{\vct{t}}_0 \quad  {\rm on} \ \Gamma_{\idx{S;N}}
\end{align}
need to be satisfied. 

For the balance equation \eqref{equ:str:momRKloc} initial conditions for displacements $\vct{d}$ and velocities 
$\dot{\vct{d}}$ need to be specified at time $t=0$,
\begin{align}
  \label{equ:str:initialbc_u}
  \vct{d}_0 \,=\, \vct{d}(\vct{X},t=0) \,=\, \hat{\vct{d}}_0 \quad {\rm on} \  \Omega_{\idx{S}} \,,
  \\ 
  \label{equ:str:initialbc_v}
  \dot{\vct{d}}_0 \,=\, \dot{\vct{d}}(\vct{X},t=0) \,=\, \hat{\dot{\vct{d}}}_0 \quad {\rm on} \  \Omega_{\idx{S}} \,,
\end{align}
where $\vct{X}$ defines the initial position.

\subsection{Fluid-structure interface conditions}
\label{sec:interfacecond}
Dynamic and kinematic coupling conditions at the conjoined interface $\Gamma$ ensure the integrity between the 
subdomains in this partitioned coupling 
algorithm.
Assuming no mass transport across the interface, normal velocities have to match, i.e.
\begin{equation}
    \vct{u}^{\Gamma} \cdot \vct{n}^{\Gamma} = \frac{\partial
\vct{d}^{\Gamma}}{\partial t} \cdot \vct{n}^{\Gamma} \quad \text{on} \ \Gamma \, ,
    \label{eq:interface_velocity_inviscid}
\end{equation}
where $\vct{n}^{\Gamma}$ denotes the interface unit normal vector.
The dynamic condition requires the tractions to be equal,
\begin{equation}
    \tns{\sigma}^{\Gamma}_\idx{F} \cdot \vct{n}^{\Gamma} = \tns{\sigma}^{\Gamma}_\idx{S} \cdot
\vct{n}^{\Gamma} \quad \text{on} \ \Gamma \, ,
    \label{eq:interface_stress_viscous}
\end{equation}
where $\tns{\sigma}_\idx{F} = -p\,\tns{I}$ denotes the 
fluid stress tensor comprising only contributions due to the pressure in the inviscid case considered here. The Cauchy 
stress tensor $\tns{\sigma}_\idx{S}$ is
defined as
\begin{equation}
 \tns{\sigma}_\idx{S} = \frac{1}{J} \tns{P} \cdot \tns{F}^{\sf T}
\end{equation}
in which $J$ is the Jacobian.

\section{Numerical approach: Fluid}
\label{sec:numerical_approach_fluid}

We employ the FVM for solving the Euler equations on 
Cartesian grids. The time-dependent fluid-solid interface conditions on $\Gamma$ are imposed by a cut-element based 
IBM.

\subsection{Mathematical model}

\begin{figure} 
  \centering
  \includegraphics[width= 0.70\columnwidth]{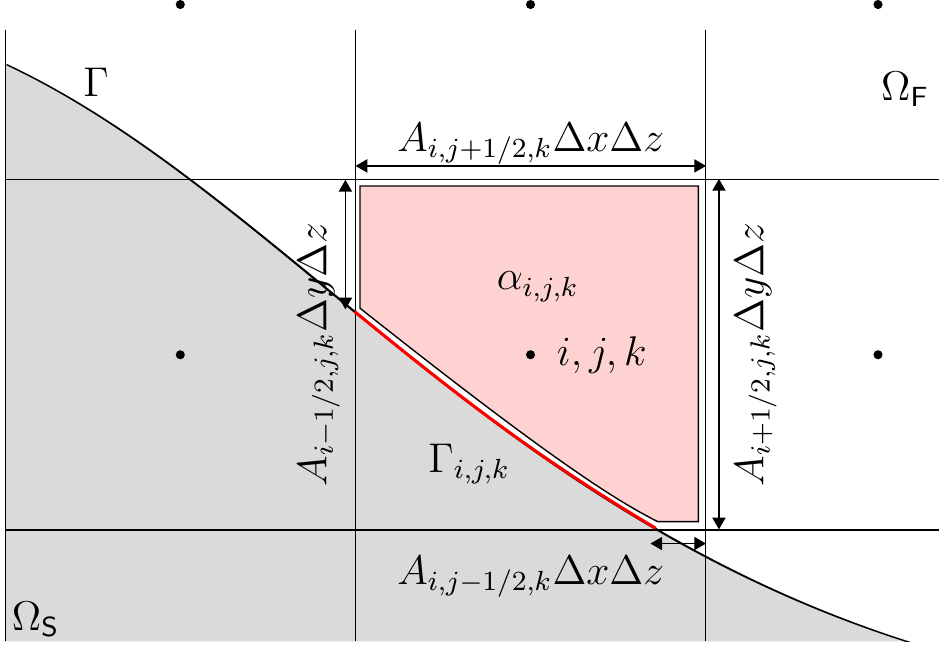}
  \caption{Two-dimensional sketch of a cut-cell $(i,j,k)$ \cite{Oerley:2015}.}
  \label{fig:ib_overview} 
\end{figure}
A sketch of a two-dimensional cut-cell is shown in 
Fig.~\ref{fig:ib_overview}. In the following, $\Gamma$ denotes the 
fluid-structure interface of the continuous problem, and $\Gamma_{\idx{F}/\idx{S}}$ the flow 
and structure side of the interface of the discrete problem. We solve the integral form of \eqref{eqn:nse},
\begin{equation}
\int_{t^n}^{t^{n+1}}\int_{\Omega_{i,j,k}\cap\Omega_{\idx{F}}}\left(\frac{\partial\vct{w}}{
\partial 
t}+\nabla\cdot\tns{K}(\vct{w})\right)\mathrm{d}x\mathrm{d}y\mathrm{d}z\ \mathrm{d}t=0,
\label{eqn:nse_int}
\end{equation}
where the integral is taken over the volume $\Omega_{i,j,k}\cap\Omega_{\idx{F}}$ of a 
computational cell $(i,j,k)$ and time step $\Delta t = t^{n+1}-t^n$. Applying the Gauss theorem results in
\begin{equation}
\int_{t^n}^{t^{n+1}}\int_{\Omega_{i,j,k}\cap\Omega_{\idx{F}}}\frac{\partial\vct{w}}{\partial 
t}\ \mathrm{d}V\ 
\mathrm{d}t+\int_{t^n}^{t^{n+1}}\int_{\partial(\Omega_{i,j,k}\cap\Omega_{\idx{F}})}\tns{K}
(\vct{w} )\cdot\vct
{n}\ \mathrm{d}S\ \mathrm{d}t=0,
\label{eqn:nse_int_gauss}
\end{equation}
where $\partial(\Omega_{i,j,k}\cap\Omega_{\idx{F}})$ denotes the wetted surface of a computational 
cell $(i,j,k)$, 
and $\mathrm{d}V$, $\mathrm{d}S$ the infinitesimal volume and surface element, respectively. Applying a volume average 
of the conserved variables 
\begin{equation}
\discret{w}_{i,j,k}=\frac{1}{\alpha_{i,j,k} 
V_{i,j,k}}\int_{\Omega_{i,j,k}\cap\Omega_{\idx{F}}}\vct{w}\ 
\mathrm{d}x\mathrm{d}y\mathrm{d}z,
\end{equation}
leads to
\begin{equation}
\begin{split}
\alpha^{n+1}_{i,j,k}\discret{w}_{i,j,k}^{n+1}&\\
=&\alpha^{n}_{i,j,k}\discret{w}_{i,j,k}^{n}\\
+&\frac{\Delta t}{ \Delta x_i} 
\left[A_{i-1/2,j,k}^n\discret{f}_{i-1/2,j,k}-A_{i+1/2,j,k}^n\discret{f}_{i+1/2,j,k}\right] \\
+&\frac{\Delta t}{ \Delta y_j} 
\left[A_{i,j-1/2,k}^n\discret{g}_{i,j-1/2,k}-A_{i,j+1/2,k}^n\discret{g}_{i,j+1/2,k} \right] \\
+&\frac{\Delta t}{ \Delta z_k} 
\left[A_{i,j,k-1/2}^n\discret{h}_{i,j,k-1/2}-A_{i,j,k+1/2}^n\discret{h}_{i,j,k+1/2} \right] \\
+&\frac{\Delta t}{ V_{i,j,k}}  \discret{\chi}_{i,j,k}.
\end{split} 
\label{eqn:vol_discr_cut}
\end{equation}
$V_{i,j,k}=\Delta x_i\Delta y_j\Delta z_k$ corresponds to the total volume of cell 
$\Omega_{i,j,k}$, $\alpha_{i,j,k}$ 
corresponds to the fluid volume fraction, $\discret{w}_{i,j,k}$ is the vector of 
volume-averaged conserved quantities in the cut-cell, and $A$ is the 
effective fluid wetted cell face aperture. The face averaged numerical fluxes across the cell faces are denoted as 
$\discret{f}, \discret{g}$ and $\discret{h}$. The 
flux $\discret{\chi}_{i,j,k}$ across the interface $\Gamma_{i,j,k}=\Gamma \cap \Omega_{i,j,k}$ is discussed in detail 
below.

Time integration of the state vector is shown here for a forward Euler time integration scheme 
with a time step $\Delta t$, which corresponds to one sub-step of an explicit Runge-Kutta method.
Appropriate initial and boundary conditions are prescribed on the domain $\Omega_{\idx{F}}$ and the surface  
$\partial\Omega_{\idx{F}}$. For all simulations presented in this paper we employ a spatial flux discretization on 
local characteristics by an $5$th-order WENO scheme \cite{Liu1994} together with a Lax-Friedrichs 
flux function. A $3$rd-order strongly stable Runge-Kutta scheme \cite{Gottlieb:1998we} 
is used for time integration.

\subsection{Conservative immersed boundary method}
 
\subsubsection{Geometry computation}
\label{sec:geometry_computation}
\begin{figure}
  \centering
  \includegraphics[width=0.9\columnwidth]{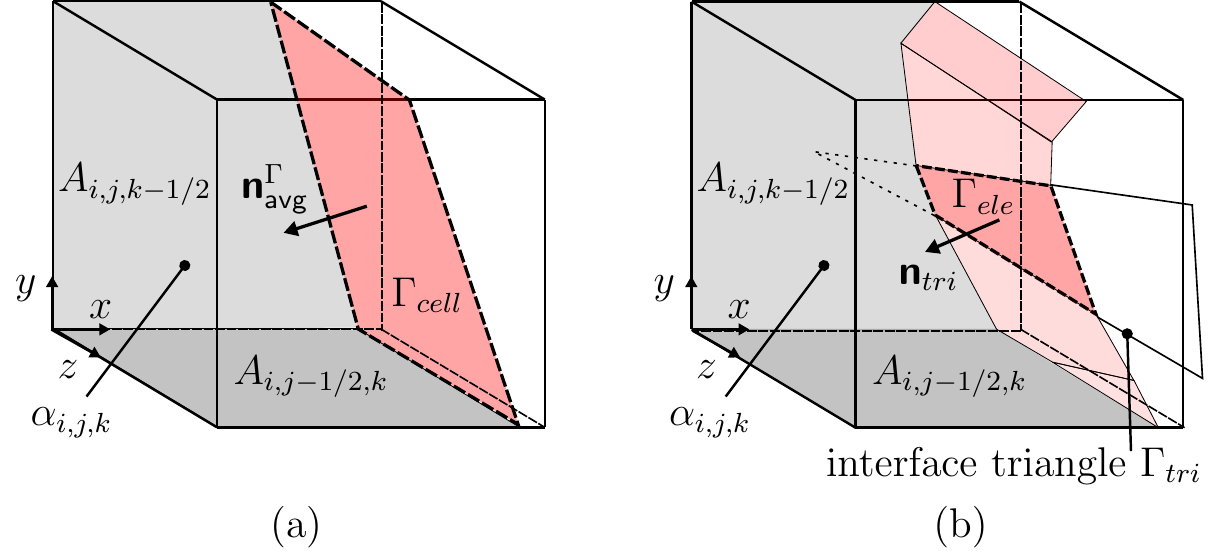}
  \caption{Computation of cut-cell properties based on a level-set field $\Phi$ (a) and on
intersection with a provided surface triangulation (b). For a detailed 
description of the cut algorithm please refer to \cite{Oerley:2015}.}
  \label{fig:cutcell}
\end{figure}

\begin{figure}
  \centering
  \includegraphics[width= 1.0\columnwidth]{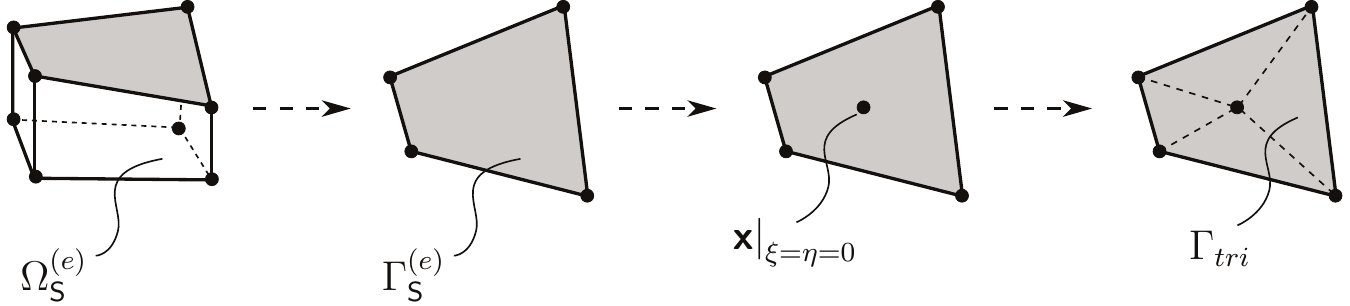}
  \caption{Triangulation of an eight-node hexahedral element face (gray) contributing to the fluid-structure 
interface $\Gamma$.}
  \label{fig:quad8_triangulation}
\end{figure}
Moving boundaries with sharp corners and complex geometries may cause numerical 
artifacts in terms of spurious pressure oscillations. Following \citet{Oerley:2015}, these 
artifacts are mainly caused by a discontinuous evolution of fluid volume fractions when utilizing a level-set based 
representation of the interface.
A solution to overcome these problems is to use an accurate representation of the geometry based 
on the computational fluid mesh and the provided structural interface. 
While the level-set method results in a planar approximation of the interface segment $\Gamma_{cell}$ in a cell, see 
Fig.~\ref{fig:cutcell}(a), the cut-element method recovers sub-cell interface resolution by a set of cut-elements 
$\Gamma_{ele}$ in a single fluid cell, see Fig.~\ref{fig:cutcell}(b). The computation of the fluid volume 
fraction $\alpha_{i,j,k}$ is done by a sub-tetrahedralization of the fluid volume, while face apertures such as 
$A_{i,j-1/2,k}$ are calculated using a sub-triangulation of the cell faces \cite{Oerley:2015}.

A linear approximation of the possibly nonlinear structural interface is used for the 
cut algorithm as an input. The element surface $\Gamma_{\idx{S}}^{(e)}$ of an 
eight-node linear brick element, which contributes to the fluid-structure interface, is highlighted in 
gray, see Fig.~\ref{fig:quad8_triangulation}. This surface is split 
into four interface triangles $\Gamma_{tri}$ using an
additional node at $\left.\discret{x}\right |_{\xi = \eta = 0}$ for improved approximation 
of its bilinear shape.

\subsubsection{Interface exchange term}
\label{sec:exchange_term}

Interaction of the fluid with a solid interface is modeled by an interface exchange term 
$\vct{\chi}_{i,j,k}$, as 
introduced in Eq.~\eqref{eqn:vol_discr_cut}. Following the approach introduced above, we can write the 
interface exchange term as a sum of all individual contributions of all cut-elements contained within this 
computational cell,
\begin{equation}
\vct{\chi}_{i,j,k}=\sum_{ele}\vct{\chi}_{ele}.
\label{eqn:X_total}
\end{equation}
For inviscid flows, the cut-element based interface exchange term $\vct{\chi}_{ele}$ 
accounts for the pressure and 
pressure work at the element interface
\begin{equation}
\vct{\chi}_{ele}=\left[
\begin{array}{c}
0 \\
\discrets{p}^{\Gamma}_{ele}\Delta \Gamma_{ele}\ \discrets{n}_1^{\Gamma;ele}\\
\discrets{p}^{\Gamma}_{ele}\Delta \Gamma_{ele}\ \discrets{n}_2^{\Gamma;ele}\\
\discrets{p}^{\Gamma}_{ele}\Delta \Gamma_{ele}\ \discrets{n}_3^{\Gamma;ele}\\
\discrets{p}^{\Gamma}_{ele}\Delta \Gamma_{ele}\ 
\left(\discret{n}^{\Gamma;ele}\cdot\discret{u}^{\Gamma;ele}\right)
\end{array}
\right]\, ,
\label{eq:ib_pressure_term}
\end{equation} 
where $\Delta \Gamma_{ele}$ is the element interface area, 
$\discret{n}^{\Gamma;ele}=[\discrets{n}_1^{\Gamma;ele},\discrets{n}_2^{\Gamma;ele},\discrets{n}_3^{\Gamma;ele}]$ is the 
element unit normal vector obtained directly from the structural interface triangle $\Gamma_{tri}$, and 
$\discret{u}^{\Gamma;ele}$ is the interface velocity evaluated at the cut-element face centroid. The element 
interface pressure $\discrets{p}^{\Gamma}_{ele}$ is obtained by solving a 
symmetric face-normal Riemann problem
\begin{equation}
 \mathcal{R}\left(\discret{w}_{i,j,k},\discret{u}^{\Gamma;ele} \right) = 0
 \label{eq:riemann_problem}
\end{equation}
for each cut-element within the cut-cell $\left(i,j,k \right)$. The exact solution of the 
reflective boundary Riemann problem \eqref{eq:riemann_problem} consists of either two shock waves 
($\discret{u}_{i,j,k}\cdot\discret{n}^{\Gamma;ele} < \discret{u}^{\Gamma;ele}\cdot\discret{n}^{\Gamma;ele}$) or two 
rarefaction waves ($\discret{u}_{i,j,k}\cdot\discret{n}^{\Gamma;ele} \geq 
\discret{u}^{\Gamma;ele}\cdot\discret{n}^{\Gamma;ele}$), which are symmetric about the path of the moving interface 
coinciding with the contact wave \cite{Toro2009}. 
The exact solution for the interface pressure $\discrets{p}^{\Gamma}_{ele}$
is the root of
\begin{equation}
 \left(\discrets{p}^{\Gamma}_{ele} - \discrets{p}_{i,j,k} \right) \cdot 
\sqrt{\dfrac{\frac{2}{\left(\gamma +1\right)\discrets{\rho}_{i,j,k}}}{\discrets{p}^{\Gamma}_{ele} + \frac{\gamma 
-1}{\gamma +1}\discrets{p}_{i,j,k}}} + \left(\discret{u}_{i,j,k}\cdot\discret{n}^{\Gamma;ele} -  
\discret{u}^{\Gamma;ele}\cdot\discret{n}^{\Gamma;ele}\right) = 0
 \label{eq:two_shock}
\end{equation}
for the two-shocks configuration, and
\begin{equation}
 \discrets{p}^{\Gamma}_{ele} = \discrets{p}_{i,j,k} \cdot \left[ 1 + 
\left(\discret{u}^{\Gamma;ele}\cdot\discret{n}^{\Gamma;ele} - 
\discret{u}_{i,j,k}\cdot\discret{n}^{\Gamma;ele} \right)\cdot \frac{\gamma-1}{2\sqrt{\gamma \discrets{p}_{i,j,k} 
/ \discrets{\rho}_{i,j,k} }} \right]^{\frac{2\gamma}{\gamma-1}}
 \label{eq:two_expansion}
\end{equation}
for the two-rarefactions configuration.

\subsubsection{Boundary conditions for solid walls}
\label{sec:ghostpoint_extend}
\begin{figure}
  \centering
    \includegraphics[width= 0.6\columnwidth]{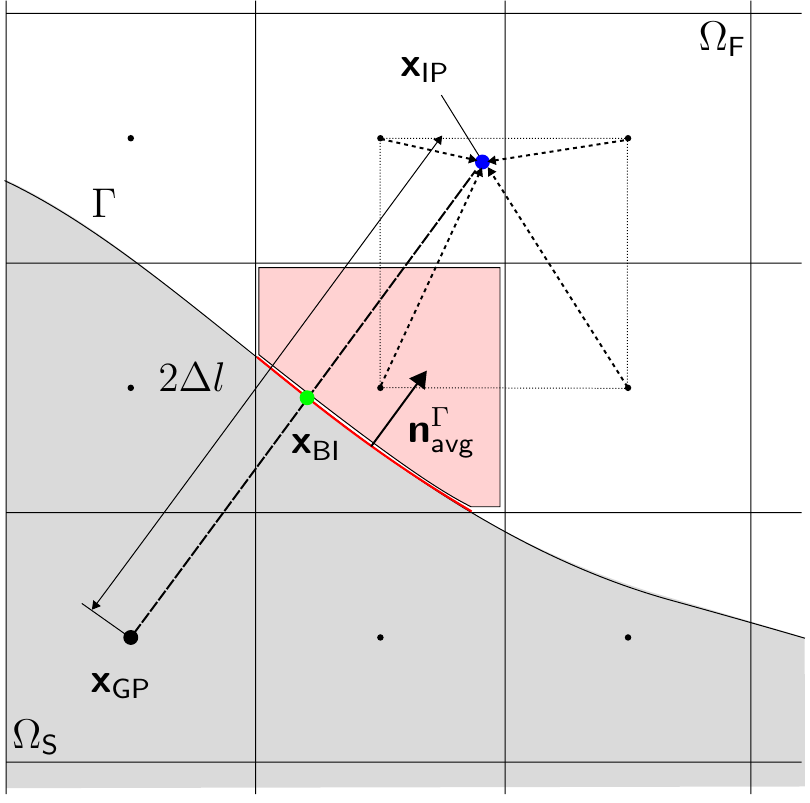}
    \caption{Construction of the ghost-cell extending procedure for a cut-cell $(i,j,k)$.}
  \label{fig:ib_extend}
\end{figure} 
Non-cut cells in the solid part of the computational domain in the vicinity of the interface contain ghost fluid states 
for imposing boundary conditions at the interface without requiring a modification of interpolation stencils in the 
finite volume reconstruction scheme. For this 
purpose, we apply the ghost-cell methodology as proposed by \citet{Mittal:2008bq}, 
extended to stationary and moving boundary cut-cell methods.
Finding the ghost-cells and extending the fluid solution across the interface does not require the fully detailed 
cut-cell geometry. We perform this procedure based on the average face centroid and normal vector of the 
cut-cell, which is an average of all contained cut-elements weighted by their area. In a first step, ghost-cells 
$\discret{x}_{\idx{GP}}$ 
that contribute to the interpolation stencil of the baseline discretization are identified, see 
Fig.~\ref{fig:ib_extend}. Next, for each ghost-cell the boundary intercept point 
$\discret{x}_{\idx{BI}}$ is computed 
such that the line segment $\overline{\discret{x}_{\idx{GP}}\discret{x}_{\idx{BI}}}$ intersects the immersed boundary 
in 
$\discret{x}_{\idx{BI}}$ normal to the interface segment. The line segment is extended into the fluid region to find 
the 
image point
\begin{equation}
\discret{x}_{\idx{IP}}=\discret{x}_{\idx{BI}}+\discret{n}^{\Gamma}_\idx{avg} \cdot \Delta l,
\end{equation}
where $\Delta l=||\discret{x}_{\idx{BI}}-\discret{x}_{\idx{GP}}||$ denotes the distance between the ghost-cell and the 
boundary intercept. Once the image point has been identified, a bilinear (in 2-D) or trilinear (in 3-D) interpolation 
is used for calculating the value of a quantity $\varphi_{\idx{IP}}$ at the image point $\discret{x}_{\idx{IP}}$:
\begin{align}
 \text{3D}:\quad \varphi(x^{\star}, y^{\star}, z^{\star}) &= c_1 + c_2 x^{\star} + c_3 y^{\star} + c_4 z^{\star} + c_5 
x^{\star}y^{\star} + c_6 x^{\star}z^{\star} \notag \\
                                          &\quad+ c_7 y^{\star}z^{\star} + c_8 x^{\star} y^{\star} z^{\star} \notag \\
 \text{2D}:\qquad\:\, \varphi(x^{\star}, y^{\star}) &= c_1 + c_2 x^{\star} + c_3 y^{\star} + c_4 x^{\star}y^{\star} \, ,
 \label{eq:image_point_interpol}
\end{align}
where $\vct{x}^{\star} = \vct{x} - \discret{x}_{\idx{IP}}$ is the relative distance vector and $\vct{c} = 
\{c_i\}$ 
are the unknown coefficients. As shown in Fig.~\ref{fig:ib_extend}, the four (eight in 3-D) coefficients can be 
determined from the variable values of the four (eight in 3-D) surrounding neighboring points,
\begin{equation}
 \vct{c} = \tns{V}^{-1} \discret{\varphi}\, ,
 \label{eq:image_point_solution}
\end{equation}
where $\discret{\varphi}$ denotes the solution at regular fluid data points and $\tns{V}^{-1}$ the inverse Vandermonde 
matrix, which is calculated by LU decomposition. After solving for \eqref{eq:image_point_solution}, the value at 
the image point is given by
\begin{equation}
 \discrets{\varphi}_{\idx{IP}} = c_1 + \mathcal{O}(\Delta^2)\, .
 \label{eq:ip_value}
\end{equation}
Ghost-cell values are obtained using a linear approximation along the 
line $\overline{\discret{x}_{\idx{GP}}\discret{x}_{\idx{BI}}}$ 
that satisfies the boundary conditions at the boundary intercept location $\discret{x}_{\idx{BI}}$. For 
Dirichlet boundary conditions, ghost-cell data are obtained as
\begin{equation}
 \varphi_{\idx{GP}} = 2\cdot \varphi_{\idx{BI}} - \varphi_{\idx{IP}} + \mathcal{O}(\Delta l^2)\, ,
 \label{eq:DBC}
\end{equation}
whereas Neumann boundary conditions are imposed as
\begin{equation}
 \varphi_{\idx{GP}} = \varphi_{\idx{IP}} - 2\cdot \Delta l\left.\left(\nabla \varphi 
\cdot \discret{n}^{\Gamma}_\idx{avg}\right)\right|_{\discret{x}_{\idx{BI}}} + 
\mathcal{O}(\Delta l^2)\, .
 \label{eq:NBC}
\end{equation}

The $5$th-order WENO scheme used in this paper requires at least 
three layers of ghost-cells to be filled. This, in turn, poses a limitation of the current framework to structures with 
a 
size larger than several fluid cells in order to fill the ghost-cell values properly. An adaptive mesh refinement 
procedure for the flow solver or the decoupling of the ghost-cell method from the underlying Cartesian grid could 
resolve this limitation.

\subsubsection{Treatment of small cut-cells} 
The time step $\Delta t$ is adjusted dynamically according to the CFL condition based on full cells of 
the underlying Cartesian grid. A drawback of cut-cell methods is that the fluid volume fraction of cut-cells may become 
very small and therefore can lead to numerical instability or require excessively small time steps with explicit time 
integration schemes and poor 
convergence with implicit methods. A stabilization of the underlying scheme is therefore required.
We employ a so-called mixing procedure as proposed in \cite{Hu2006, Oerley:2015}.

\section{Numerical approach: Solid}
\label{sec:str:fem}
The FEM is applied to solve the structural problem. Hence, we start with the weak form of the
structural field equation, which is obtained by building weighted residuals of the balance
equation~\eqref{equ:str:momRKloc} and Neumann boundary conditions~\eqref{equ:str:nbc} with virtual displacements
$\delta \vct{d}$. Subsequently, the divergence theorem is applied, yielding 
\begin{equation}\label{equ:weimp}
\begin{array}{ll}
  \displaystyle \int_{\Omega_{\idx{S}}} \rho_{\idx{S};0} \, \ddot{\vct{d}} \, \cdot \delta \vct{d} \, \mathrm{d}V_0 
  \, + \, \displaystyle \int_{\Omega_{\idx{S}}} \tns{S} : \delta \tns{E} \, \mathrm{d}V_0 \,-\, \\[2ex]
  \quad -\, \displaystyle \int_{\Omega_{\idx{S}}} \hat{\vct{b}}_0 \cdot \delta \vct{d} \, \mathrm{d}V_0
  \,-\, \int_{\Gamma_{{\rm N};\idx{S}}} \, \hat{\vct{t}}_0 \cdot \delta \vct{d} \, \mathrm{d}A_0
  \,-\, \delta \mathcal{W}^\Gamma_\idx{S}
  \,=\, \vct{0}
\end{array}
\end{equation}
with infinitesimal volume and surface elements, $\mathrm{d}V_0$ and $\mathrm{d}A_0$, respectively.
Herein, $\delta \tns{E}$ is obtained as result of the variation of the Green-Lagrange
strain~\eqref{equ:GLstrain}, i.e.
\begin{equation}
  \delta\tns{E} \,=\, \dfrac{1}{2} \, \Big( (\nabla_0 \, \delta \vct{d} )^{\sf T} \cdot \tns{F} 
  \,+\,
  \tns{F}^{\sf T} \cdot
  \nabla_0 \, \delta {\vct{d}} \Big )
\end{equation}
with $\nabla_0 \,(\bullet)$ representing the material gradient operator. The influence of the interface on the
structure is introduced via the additional virtual work term $\delta \mathcal{W}^\Gamma_\idx{S}$.

The weak form of equation~\eqref{equ:weimp} is discretized in space with the FEM.
The solid domain $\Omega_{\idx{S}}$ is split into $n^\idx{e}$ elements $\Omega_{\idx{S}}^{(e)}$ (subdomains). The 
semi-discrete weak form of the balance of linear momentum is obtained by assembling the contributions of all elements, 
leading to
\begin{equation}\label{equ:str:semidisweakvirt}
  \discret{M} \, \ddot{\discret{d}} 
  \,+\, \discret{f}_{\idx{S;int}}(\discret{d})
  \,-\, \discret{f}_{\idx{S;ext}}(\discret{d})
  \,-\, \discret{f}^\Gamma_\idx{S}
  \,=\, \discret{0} \,,
\end{equation}
where we have assumed the discrete virtual displacement vector $\delta \discret{d}$ to be arbitrary. The
vectors $\ddot{\discret{d}}$ and $\discret{d}$ describe the discrete acceleration and
displacement vectors, respectively, $\discret{M}$ denotes the mass matrix, $\discret{f}_{\idx{S;int}}$ and
$\discret{f}_{\idx{S;ext}}$ the internal and external force vectors. The interface traction of the fluid on the
structure is described by $\discret{f}^\Gamma_\idx{S}$. Element technology such as the method of enhanced assumed 
strains (EAS), as introduced in \cite{SimoAr1993}, is used in order to avoid locking phenomena. For time 
integration, the generalized trapezoidal rule (or one-step-$\theta$ scheme) is employed for the structure solver in this 
work. Thus, applying this scheme to the
semi-discrete equation~\eqref{equ:str:semidisweakvirt}, the final fully discrete structural equation at the new time
step $n+1$ is obtained.

The fully discrete structural equation describes a system of nonlinear algebraic equations which is solved iteratively 
by a Newton-Raphson method. The linearized system reads
\begin{equation}\label{equ:str:linsyst}
  \discret{K}_{\idx{SS}}(\discret{d}^{n+1}_i) \, \Delta \discret{d}_{i+1}^{n+1}
  \,=\, -\discret{r}_{\idx{S}}(\discret{d}_i^{n+1}) 
\end{equation}
with iteration step $i$, the dynamic effective structural stiffness matrix $\discret{K}_{\idx{SS}}$, and the 
residual vector $\discret{r}_{\idx{S}}$. Thus, a new solution of the displacement increment $\Delta 
\discret{d}_{i+1}^{n+1}$ for current iteration step $i+1$ is determined, and the final displacement solution of time 
step $n+1$ is obtained via updating
\begin{equation}\label{equ:str:iterinc}
  \discret{d}_{i+1}^{n+1} \,=\, \discret{d}_i^{n+1} \,+\, \Delta \discret{d}_{i+1}^{n+1} \,.
\end{equation}
The Newton-Raphson iteration is considered as converged if $\lvert \discret{r}_{\idx{S}} \rvert_2 \leq \epsilon$ is
satisfied using a sufficiently small tolerance $\epsilon$.

\section{Coupling procedure}
\label{sec:coupling_procedure}

\subsection{Treatment of non-matching interfaces}
\label{sec:non_matching_interface}
The reconstruction of the interface on the fluid side based on the structural position 
leads to a change in the number of cut-elements in each coupling step and to a change in connectivity, which 
inevitably results in a non-matching interface. 
A Mortar method has been chosen in this work as it preserves linear and angular momentum.
The Mortar method requires the choice of a so-called slave and master side of the interface $\Gamma^\idx{sl}$ and
$\Gamma^\idx{ma}$, respectively. Primary coupling variables, such as velocities in our case, are transferred from
the master to the slave side, and secondary variables, such as tractions, are transferred vice versa. The 
Dirichlet-Neumann partitioning chosen here determines the fluid to be the slave side ($\Gamma^\idx{sl} \equiv
\Gamma_\idx{F}$) and the solid to be the master 
side ($\Gamma^\idx{ma} \equiv \Gamma_\idx{S}$) with respect to Mortar coupling. The aim is to obtain discrete projection
operators for consistent data transferring.

In the following derivation, a no-slip condition between fluid and solid is assumed instead of the slip condition in
\eqref{eq:interface_velocity_inviscid} for simplicity, which will later be 
released again. The starting point is 
the weak form of the continuity constraint 
\begin{equation} \label{equ:mortar_weak_constr}
 \delta W_\idx{\lambda} = \int_{\Gamma^\idx{sl}}{\delta \transpose{\tns{\lambda}} \left( \tns{u}^\idx{\Gamma} - 
  \dot{\tns{d}}^\idx{\Gamma} \right) \text{d} \Gamma} = 0
\end{equation}
together with weak form of the equilibrium of tractions at the interface
\begin{equation} \label{equ:mortar_virtual_work}
 \delta W_\idx{\Gamma} = \int_{\Gamma^\idx{sl}}{\transpose{\tns{\lambda}} \left( \delta \tns{u}^\idx{\Gamma} - \delta
  \dot{\tns{d}}^\idx{\Gamma} \right) \text{d} \Gamma}
\end{equation}
in which a Lagrange multiplier field $\tns{\lambda} = \tns{\sigma}^{\Gamma}_{\idx{F}}\cdot \tns{n}^{\Gamma}$
and the corresponding test functions $\delta \tns{\lambda}$ are introduced. The virtual
work term \eqref{equ:mortar_virtual_work} is the
conjugate term of \eqref{equ:mortar_weak_constr} and it contains virtual work contributions of interface
tractions on the fluid side and on the solid side, $\delta \mathcal{W}^{\Gamma}_F$ and $\delta \mathcal{W}^{\Gamma}_S$,
respectively. Additionally, $\delta \mathcal{W}^{\Gamma}_S$ needs to be adapted to the chosen time integration
scheme for the solid due to the occurrence of the time derivative of the displacements.

\begin{figure}
  \centering
  \includegraphics[width= 0.7\columnwidth]{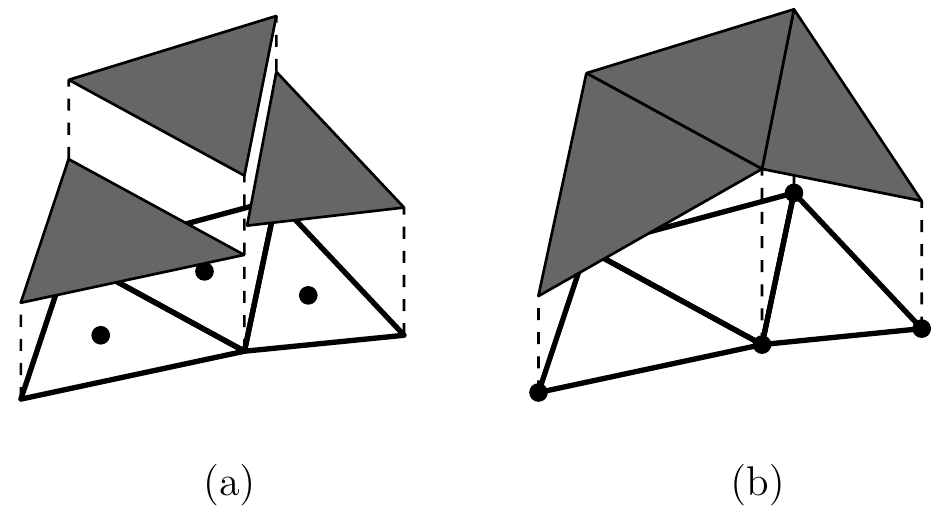}
  \caption{Interpolation of state variables. (a) FVM: constant value per cell, (b) FEM: linear Lagrange polynomials }
  \label{fig:interpol_on_struct_side}
\end{figure}

An important question is which ansatz
functions should be used for a proper interpolation of the respective fields at the interface. Due to the applied cut
procedure in the underlying finite volume discretization it is not possible to obtain the surface ansatz functions for
the cut-elements based on a trace space relationship. Without invoking high-order reconstruction, the FVM defines for 
the state values in the cut-cells a piecewise constant field as it is depicted in 
Fig.~\ref{fig:interpol_on_struct_side}(a). For the
solid, it is possible to obtain the ansatz functions from the trace space of the underlying volume element leading to an
interpolation with standard Lagrange polynomials as it is shown in Fig.~\ref{fig:interpol_on_struct_side}(b).
Hence, a Lagrange multiplier 
\begin{equation} \label{equ:lagrange_mult}
 \discret{\lambda} = \sum_{j=1}^{n^\idx{sl}}{\Phi_j \lambda_j}
\end{equation}
using constant ansatz functions $\Phi_j$ on each cut-element can be
utilized, which is defined only on the slave side of the interface. The discrete Lagrange multipliers are denoted
as $\lambda_j$. Due to the constant value in a single cut-element also the velocities can be represented using
constant ansatz functions $N^\idx{sl}_k$. This approach then reads
\begin{equation} \label{equ:mortar_vel_interpol_sl}
  \discret{u} = \sum_{k=1}^{n^\idx{sl}}{N^\idx{sl}_k \discret{u}_k} .
\end{equation}
In \eqref{equ:lagrange_mult} and \eqref{equ:mortar_vel_interpol_sl}, the total number of cut-elements is denoted with
$n^\idx{sl}$, which is equal to the number of discrete fluid
velocities $\discret{u}_k$ due to the piecewise constant field on each cut-element.
In contrast, standard shape functions $N^\idx{ma}_l$ based on Lagrange
polynomials are used for the interpolation of the velocities on the solid side of the interface. This leads
to
\begin{equation} \label{equ:mortar_vel_interpol_ma}
 \dot{\discret{d}} = \sum_{l=1}^{n^\idx{ma}}{N^\idx{ma}_l \dot{\discret{d}}_l}
\end{equation}
where the total number of discrete solid velocities
$\dot{\discret{d}}_l$ is denoted as $n^\idx{ma}$, which is equal to the number of nodes in the solid interface.
Inserting \eqref{equ:lagrange_mult} $-$~\eqref{equ:mortar_vel_interpol_ma} into~\eqref{equ:mortar_weak_constr}
leads to
\begin{equation} \label{equ:mortar_discret_constraint}
 \begin{array}{ll}
\displaystyle
  \hspace{-2cm}\delta W_\idx{\lambda} = \sum_{j=1}^{n^\idx{sl}} { \sum_{k=1}^{n^\idx{sl}} {
   \delta \discret{\lambda}_j^T \left( \int_{\Gamma^\idx{sl}} {\Phi_j N^\idx{sl}_k \, {\rm d} \Gamma} \right) 
   \discret{u}_k } }
   - \\[2ex]
\displaystyle
   \hspace{2cm} - \sum_{j=1}^{n^\idx{sl}} { \sum_{l=1}^{n^\idx{ma}} {
   \delta \discret{\lambda}_j^T \left( \int_{\Gamma^\idx{sl}} {\Phi_j N^\idx{ma}_l \, {\rm d} \Gamma} \right)
  \dot{\discret{d}}_l } } = 0 .
 \end{array}
\end{equation}
Therein, nodal blocks of the two Mortar integral matrices commonly
denoted as $\discret{D}$ and $\discret{M}$ can be identified. This leads to the following definitions:
\begin{align} \label{equ:mortar_D}
 \discret{D}[j,k] &= D_{jk} \discret{I}_3 = \int_{\Gamma^{sl}}{\Phi_j N^\idx{sl}_k \, {\rm d} \Gamma} \discret{I}_3 \, ,
 \\
 \label{equ:mortar_M}
 \discret{M}[j,l] &= M_{jl} \discret{I}_3 = \int_{\Gamma^\idx{sl}}{\Phi_j N^\idx{ma}_l \, {\rm d} \Gamma}
\discret{I}_3
\end{align}
with the $3 \times 3$ identity tensor $\discret{I}_3$, whose size is determined by the number of variables to be coupled
for each node. Here, $\discret{D}$ is a square $3\,n^\idx{sl} \times 3\,n^\idx{sl}$ 
matrix, which has only diagonal entries due to the choice of piecewise constant shape functions, whereas the
definition of $\discret{M}$ generally gives a rectangular matrix of dimensions $3\,n^\idx{sl} \times 3\,n^\idx{ma}$.
The actual numerical integration of the Mortar integrals can be performed either segment-based or element-based, 
see~\cite{Farah2014, Fischer2006, Popp2010, Puso2004b}. Due
to its superior numerical efficiency, element-based integration is used exclusively in this work. 

Plugging the previously defined Mortar matrices $\discret{D}$ and $\discret{M}$ into \eqref{equ:mortar_weak_constr}
leads to the discrete continuity constraint
\begin{equation} \label{equ:mortar_discrete_constraint}
 \discret{D} \cdot \discret{u} - \discret{M} \cdot \dot{\discret{d}} = \discret{0} ,
\end{equation}
which will be utilized in Section~\ref{sec:transfer_solid_fluid} for the specific transfer of velocities from the solid
to the fluid interface.
Similarly, inserting \eqref{equ:lagrange_mult} $-$~\eqref{equ:mortar_vel_interpol_ma}
into~\eqref{equ:mortar_virtual_work} and again using \eqref{equ:mortar_D} and \eqref{equ:mortar_M} results in
\begin{align} \label{equ:f_gamma_F}
 \discret{f}^\Gamma_\idx{F} &= \transpose{\discret{D}} \discret{\lambda} \, ,
 \\
 \label{equ:f_gamma_S}
 \discret{f}^\Gamma_\idx{S} &= \transpose{\discret{M}} \discret{\lambda} \, ,
\end{align}
which defines the nodal coupling forces $\discret{f}^\Gamma_\idx{F}$ and $\discret{f}^\Gamma_\idx{S}$ of
the fluid and the solid, respectively. The transfer of loads is based on \eqref{equ:f_gamma_F} and
\eqref{equ:f_gamma_S} and will be described in Section~\ref{sec:transfer_fluid_solid}.

\subsubsection{Transfer of solid velocities to fluid interface}
\label{sec:transfer_solid_fluid}
 The velocity at the cut-element face centroid is needed for both the energy equation and for determining the 
 interface pressure $\discrets{p}^{\Gamma}_{ele}$ through a Riemann solver, see \eqref{eq:ib_pressure_term}. Moreover,
the kinematic constraint  \eqref{eq:interface_velocity_inviscid} requires matching normal velocities at the interface.
In a first step, the full interface velocities are transferred to the fluid by reordering
\eqref{equ:mortar_discrete_constraint} and defining a discrete projection $\discret{P}$ operator, viz. 
\begin{equation} \label{equ:discret_proj_operator_SF}
 \discret{u} = \discret{D}^{-1} \cdot \discret{M} \ \dot{\discret{d}} = \discret{P} \ \dot{\discret{d}}.
\end{equation}
 It shall be noted that the inversion of $\discret{D}$ is a trivial task at negligible cost due to its diagonal shape
and thus there is no need for solving a possibly large linear system.
In a second step, the current normal direction of the cut-element is used to project the velocity to fulfill the slip
condition.

\subsubsection{Transfer of fluid forces to solid interface}
\label{sec:transfer_fluid_solid}
The equilibrium of forces requires the surface tractions of fluid and solid to be equal. 
As we do not want to solve explicitly for the Lagrange multipliers we reorder \eqref{equ:f_gamma_F} and
\eqref{equ:f_gamma_S}, yielding
\begin{equation} \label{equ:mortar_proj_operator}
 \discret{f}^\Gamma_\idx{S} = \transpose{ \left( \discret{D}^{-1} \cdot \discret{M} \right) } \discret{f}^\Gamma_\idx{F}
= \transpose{\discret{P}} \discret{f}^\Gamma_\idx{F} \, .
\end{equation}
One can see that the transfer of loads from the fluid to the solid is based on the transpose of the projection 
operator for the transfer of solid velocities to the fluid. This is a crucial requirement 
for the consistent transfer across the interface and a distinctive feature of Mortar methods.

\subsection{Loosely coupled partitioned FSI algorithm}
\label{sec:loose_coupling}
\begin{figure}
  \centering
  \includegraphics[width= 0.7\columnwidth]{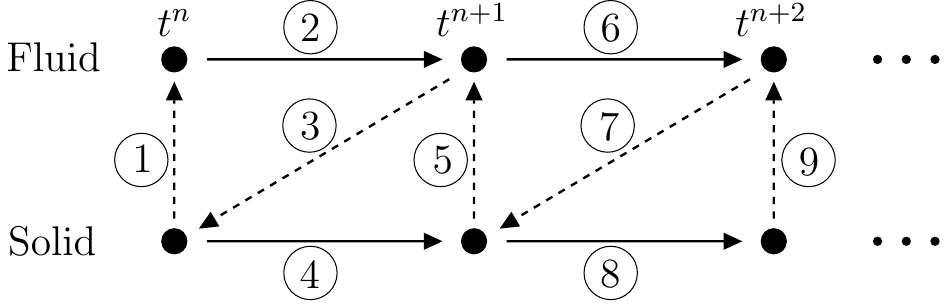}
  \caption{Schematic of the staggered time integration of the coupled system.}
  \label{fig:css_algorithm}
\end{figure}

In this paper, we use a loosely coupled conventional serial staggered algorithm. In Fig.~\ref{fig:css_algorithm}, we 
illustrate the main steps to advance the coupled system from time level $t^{n}$ to $t^{n+1} = t^{n} + \Delta t^{n}$. 
This explicit staggering algorithm, which follows the classical Dirichlet-Neumann partitioning, reads as follows:
\begin{itemize}
 \item[1.] 
 The known structural interface displacements $\discret{d}^{\Gamma;n}$ and velocities
$\dot{\discret{d}}^{\Gamma;n}$ at time $t^{n}$ are used to update the cut-cells list and geometric properties on the 
fluid side. 
For this purpose, the cut-element algorithm is applied on the triangulated structural interface (see 
Fig.~\ref{fig:quad8_triangulation}).
 \item[2.] 
 Advance the fluid in time. The evaluation of the interface exchange term \eqref{eq:ib_pressure_term} and 
the prescription of ghost-cell values \eqref{eq:DBC} and \eqref{eq:NBC} at time $t^{n+1}$ use given structural 
interface velocities $\dot{\discret{d}}^{\Gamma;n}$.
An interpolation procedure is needed to transfer solid velocities to the fluid interface, see 
Section~\ref{sec:transfer_solid_fluid}.
 \item[3.] 
 Transfer the fluid interface normal tractions $\tns{\sigma}^{\Gamma;n+1}_\idx{F} \cdot 
\discret{n}^{\Gamma;n}$ due to pressure loads to the structural solver. The staggering procedure leads to a time
shift between the stress tensor and the normal used to compute the tractions. 
An interpolation procedure is needed to transfer fluid forces to the solid 
interface, see Section~\ref{sec:transfer_fluid_solid}.
 \item[4.] Advance the structure in time while the fluid interface loads act as additional Neumann boundary
condition on the solid. 
 \item[5.] Proceed to the next time step.
\end{itemize}

Using the structural displacement $\discret{d}^{\Gamma;n}$ for the fluid solution at time $t^{n+1}$
results in a first-order in time, $\mathcal{O}(\Delta t)$, coupling scheme \cite{Farhat1995}. 
Moreover, the explicit staggering algorithm is
only conditionally stable since at time level $t^{n+1}$, the continuity condition is satisfied only for the dynamic
part ($\tns{\sigma}^{\Gamma;n+1}_{\idx{F}}\cdot \discret{n}^{\Gamma;n}$ matches
$\discret{\sigma}^{\Gamma;n+1}_{\idx{S}} \cdot \discret{n}^{\Gamma;n}$). For the kinematic part, the fluid
velocities $\discret{u}^{\Gamma;n+1}$ at $t^{n+1}$ match the structural velocities $\dot{\discret{d}}^{\Gamma;n}$ 
from the previous time step, but not the current structural velocities $\dot{\discret{d}}^{\Gamma;n+1}$. This in turn 
explains the violation of energy conservation at the interface.

\newpage
\section{Validation of the FSI algorithm}
\label{sec:validation}
In the following, we present a validation of our method for rigid and deformable structures. The solution of both 
subdomains ($\Omega_{\idx{F}}$,$\Omega_{\idx{S}}$) is advanced by the same time step which is based on the 
CFL condition for the fluid flow. For all examples, coupling is performed at every time step.
\subsection{Shock wave impact on rigid cylinder}
\label{sec:shock_cylinder}
\begin{figure}
  \centering
 \includegraphics[width= 1.0\columnwidth]{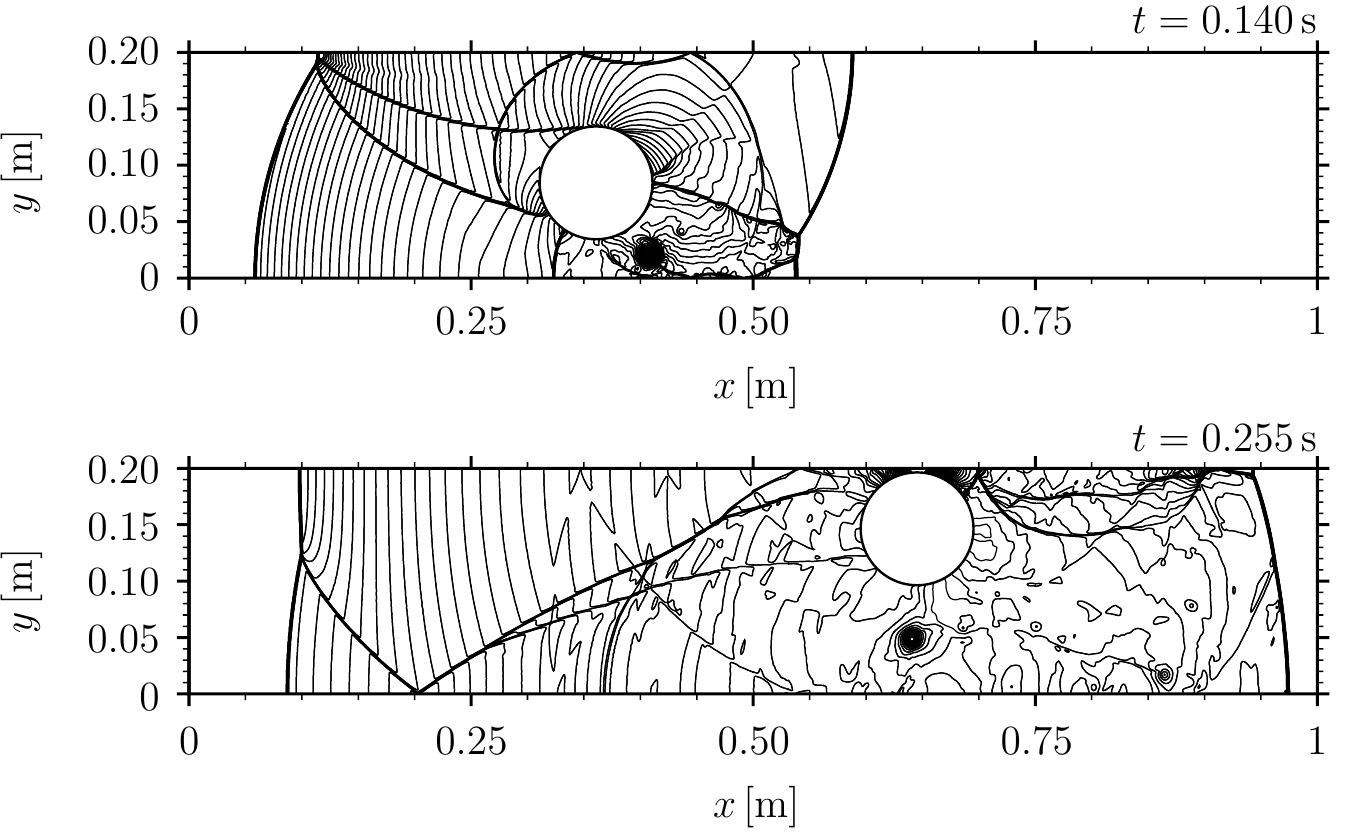}
  \caption{60 contours of fluid pressure within $0 - 28$\, Pa at two different time instances. Mesh resolution: $\Delta 
x = \Delta y = 6.25 \times 10^{-4}\,\text{m}$}
  \label{fig:shock_cylinder_pressure}  
\end{figure}
\begin{figure}[htbp]
  \centering
  \begin{subfigure}{\textwidth}
	  \centering
	  \includegraphics{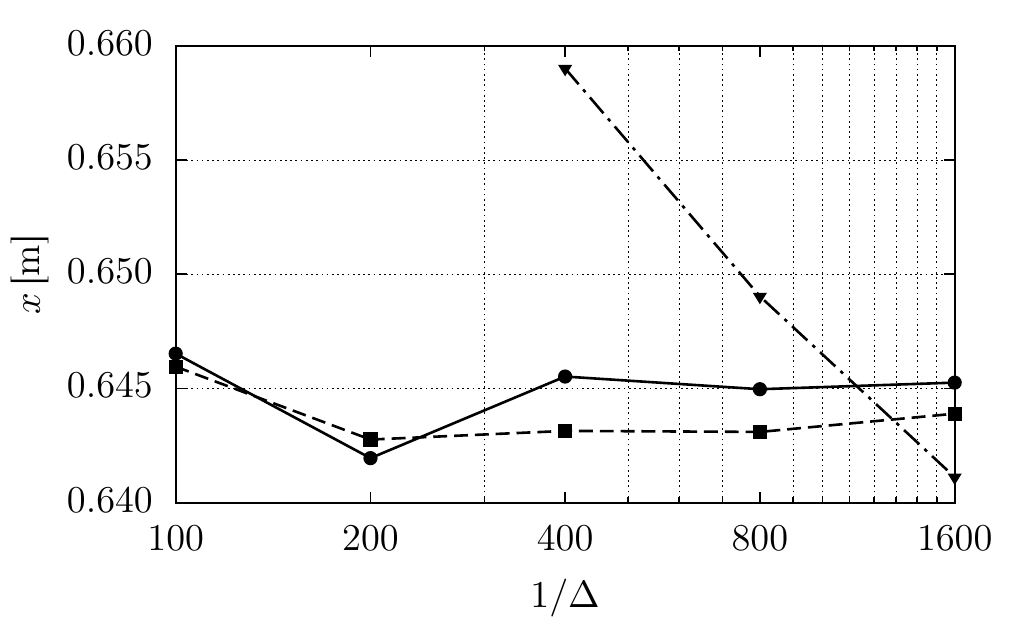}
  \end{subfigure}%
  \\
  \begin{subfigure}{\textwidth}
	  \centering
	  \includegraphics{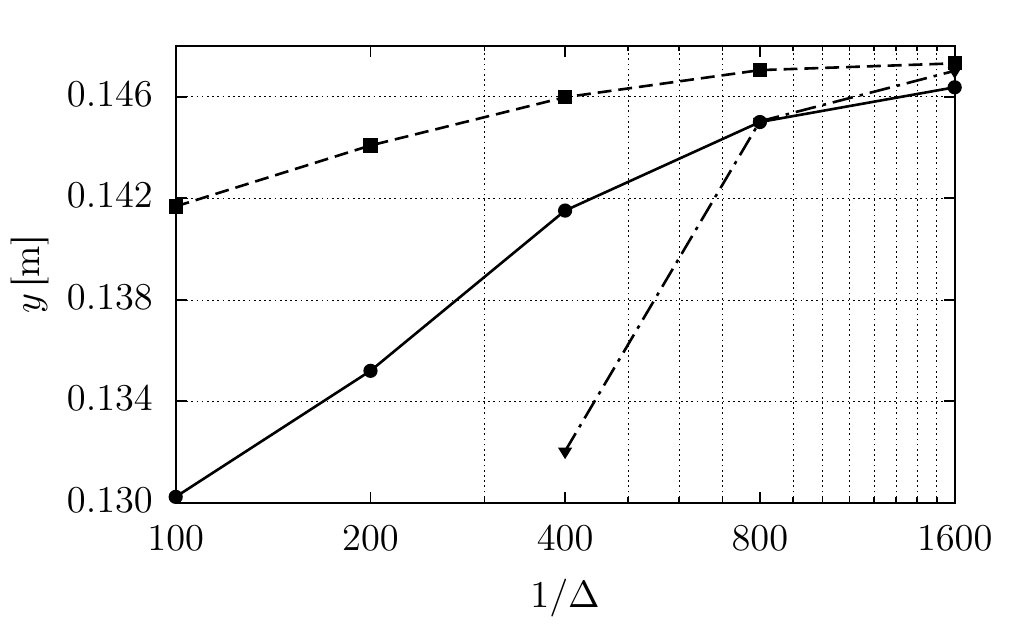}
  \end{subfigure}
  \caption{Convergence study of horizontal (top) and vertical (bottom) cylinder center position for 
different fluid mesh resolutions. (--- $\bullet$ ---) present results, ($-\cdot- \blacktriangledown -\cdot-$) 
\citet{Hu2006}, ($- -$\solidsquare$- -$) \citet{Monasse2012}. The $x$-axis is given in logarithmic scale.}
  \label{fig:shock_cylinder_convergence}
\end{figure}
The following test case for rigid body motion has been originally proposed by \citet{Falcovitz:1997wq}
and has been widely adopted in the literature, see e.g. \cite{Hu2006, Monasse2012}. The setup 
consists of a two-dimensional channel filled with air and a rigid light-weight cylinder 
of density $\rho_{\idx{S};0} = 7.6 \,\text{kg/m}^3$ initially resting on the lower wall at a position $(x,y) = 
(0.15,0.05)\,\text{m}$. The cylinder is subsequently driven and lifted upwards by a $\text{Ma} = 3$ shock wave 
entering the domain from the left. The pre-shock conditions $\rho_{\idx{F};\idx{R}} = 1\,\text{kg/m}^3,\: p_\idx{R} = 
1\, 
\text{Pa},\: u_\idx{R} = 0\, \text{m/s}$ hold for $x \ge 0.08\, \text{m}$ while for $x < 0.08\,\text{m}$ 
post-shock conditions $\rho_{\idx{F};\idx{L}} = 3.857\,\text{kg/m}^3,\: p_\idx{L} = 10.33\, \text{Pa},\: 
u_\idx{L} = 2.629\, \text{m/s}$ are initially prescribed. The fluid domain is rectangular with dimensions 
$1\,\text{m} \times 0.2\,\text{m}$ and is 
discretized with $1600 \times 320$ cells in streamwise and wall-normal direction, respectively. This leads to a grid 
resolution of $\Delta x = \Delta y = 6.25\times 10^{-4}\,\text{m}$. For the lower and upper wall, reflecting slip-wall 
boundary conditions are used. At the inflow the post-shock values are prescribed while a linear 
extrapolation of all flow variables is used at the outflow. The cylinder has a radius of $r = 
0.05\,\text{m}$ and it is discretized with $240$ tri-linearly interpolated hexahedral elements along its
circumference, leading to $240$ surface elements that are coupled to the fluid. Due to stability 
reasons the cylinder does not exactly rest on the lower wall initially. We found that a narrow gap equal to $2\%$ of 
the local cell height leads to stable and accurate results. Rigidity is achieved by imposing 
a high
Young's modulus. The time integration factor $\theta = 0.66$ is chosen for the structural time integration. A CFL number
of $0.6$ is adopted for all simulations. It should be noted that no analytical solution for the final position of the 
cylinder exists. We therefore put emphasis on convergence properties of the proposed coupling algorithm. 

\begin{figure}[htbp]
  \centering
  \includegraphics{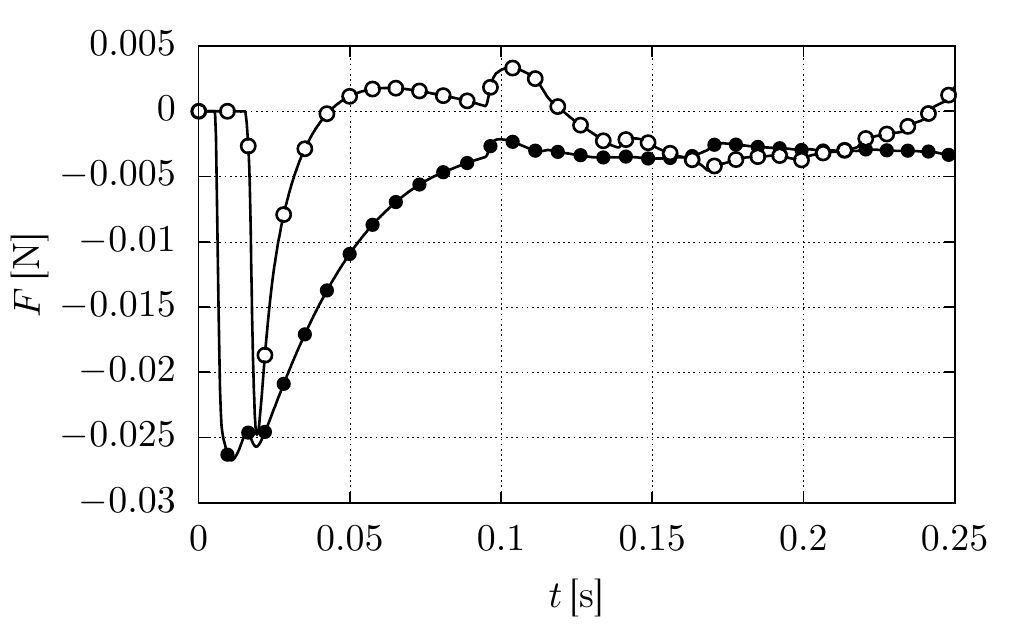}
  \caption{Temporal evolution of global forces acting on the rigid cylinder. \mbox{(--- 
$\bullet$---) $F_x$}, \mbox{(--- $\circ$ ---) $F_y$}. Mesh resolution: $\Delta x = \Delta y = 6.25 \times 
10^{-4}\,\text{m}$.}
  \label{fig:shock_cylinder_forces}
\end{figure}

Instantaneous pressure contours at $t = 0.14\,\text{s}$ and $t = 0.255\,\text{s}$ are shown in 
Fig.~\ref{fig:shock_cylinder_pressure}. With respect to the cylinder position and the resulting shock patterns our 
results agree well to Fig.~19 of \citet{Hu2006}
and Fig.~11 of \citet{Monasse2012}.  We observe a strong 
vortex beneath the 
cylinder, which persists throughout the entire cylinder trajectory, see Fig.~\ref{fig:shock_cylinder_pressure}, 
supporting the results of \cite{Forrer, Monasse2012}. By further increasing the mesh resolution up to $\Delta x = 
\Delta y = 1.5625\times 10^{-4}\,\text{m}$ the vortex is still apparent, excluding numerical dissipation being 
responsible for the formation of the vortex. As stated by \citet{Monasse2012}, a 
Kelvin-Helmholtz instability of the contact discontinuity present under the cylinder is the likely cause for this 
vortex.

Fig.~\ref{fig:shock_cylinder_convergence} shows convergence results on the final horizontal and vertical 
position of the center of mass of the cylinder together with results from literature \cite{Hu2006, Monasse2012}. The 
final position is in the same range as the results of \cite{Hu2006, Monasse2012}. Our results show a convergence 
rate similar to the results obtained by \citet{Monasse2012}.

Finally, Fig.~\ref{fig:shock_cylinder_forces} shows the temporal evolution of resulting fluid forces acting on the 
rigid cylinder, which has been obtained by summation of all individual cut-element interface exchange terms. The smooth 
force distribution confirms that our interface treatment is accurate and free of spurious pressure oscillations.

\subsection{Shock wave impact on deforming panel}
\label{sec:shock_panel}
\begin{figure}
  \centering
 \includegraphics[width= 1.0\columnwidth]{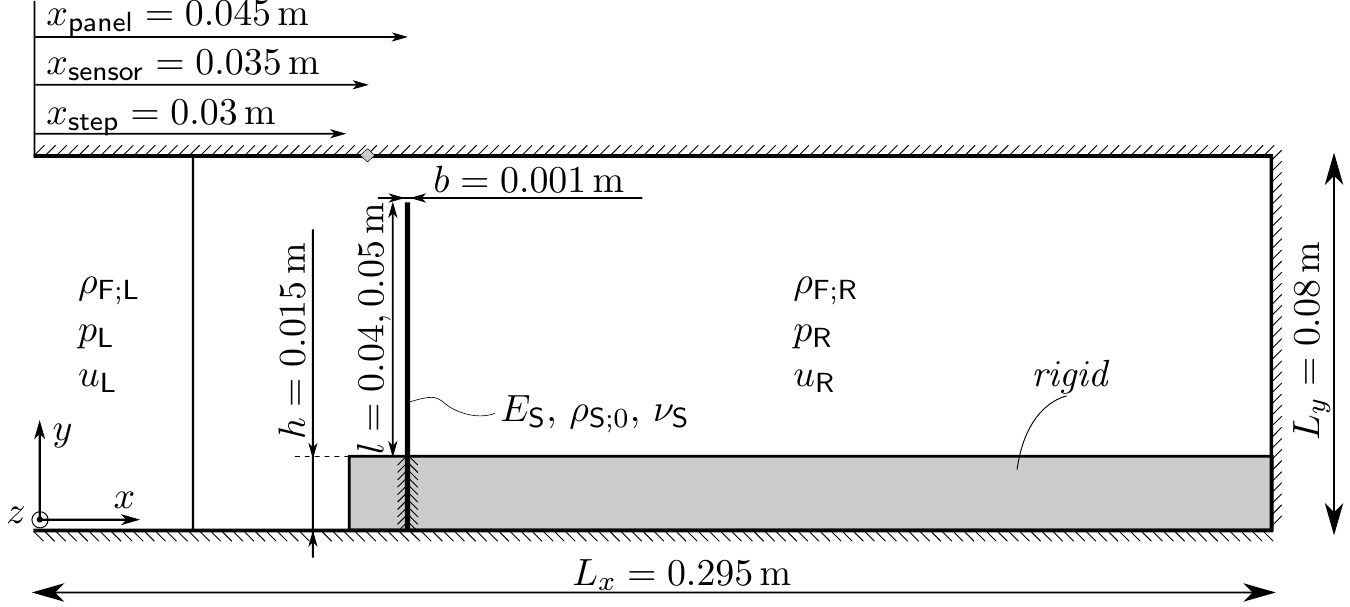}
  \caption{Setup for shock wave impact on deforming panel including geometric dimensions; see \cite{Giordano2005} for 
details.}
  \label{fig:shock_panel_setup}  
\end{figure}
\begin{figure}
  \centering
 \includegraphics[width= 1.0\columnwidth]{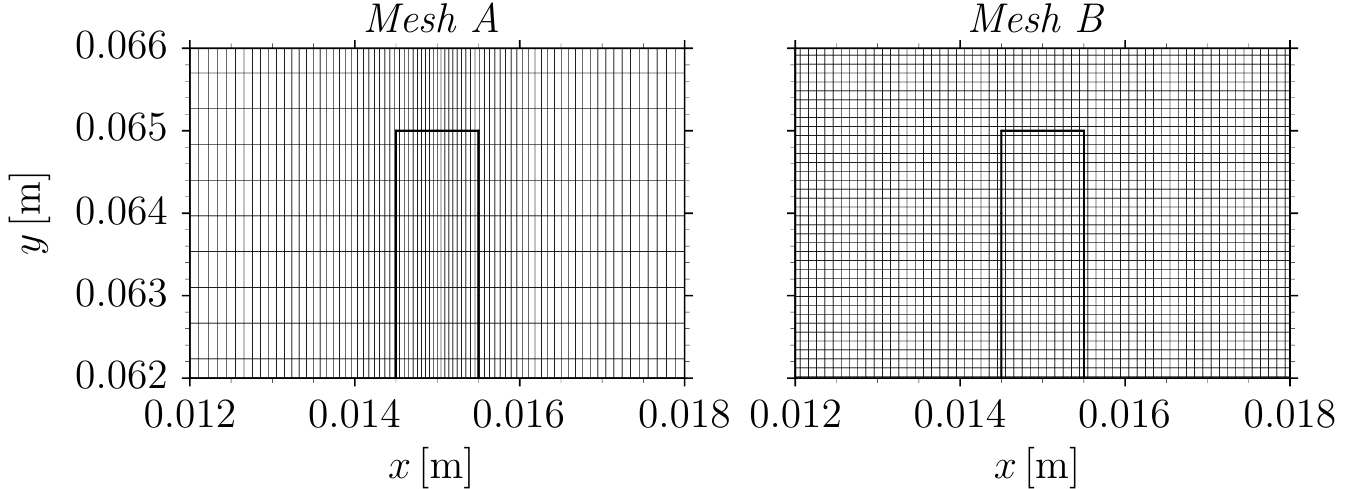}
  \caption{Fluid mesh resolutions close to the panel.}
  \label{fig:shock_panel_grid}  
\end{figure}

The behavior of a cantilever panel subjected to a shock tube flow is analyzed. This test case has been
investigated both experimentally and numerically in \cite{Giordano2005}. The experimental setup, as shown in 
Fig.~\ref{fig:shock_panel_setup}, consists 
of a deformable panel of length $l = \left[0.04,0.05\right]\,\text{m}$ and width $b = 0.001\,\text{m}$ placed within a 
shock 
tube. The panel is hit by a $\text{Ma} = 1.21$ shock wave, which enters the domain from left. The panel is made of 
steel ($E_{\idx{S}} = 220\,\text{GPa}$, $\rho_{\idx{S};0} = 7600\,\text{kg/m}^3$, $\nu_{\idx{S}} = 0.33$) and is 
clamped to a rigid 
forward-facing step at its lower end. The pre-shock conditions resemble air at rest and are set to 
$\rho_{\idx{F};\idx{R}} = 
1.189\,\text{kg/m}^3,\: p_\idx{R} = 100\, \text{kPa},\: u_\idx{R} = 0\, \text{m/s}$, while the 
post-shock values are $\rho_{\idx{F};\idx{L}} = 
1.616\,\text{kg/m}^3,\: p_\idx{L} = 154\, \text{kPa},\: u_\idx{L} = 109.68\, \text{m/s}$. The fluid 
domain is rectangular with 
dimensions $0.295\,\text{m} \times 0.08\,\text{m}$ in width and height. Since the 
problem is considered as two-dimensional, we adopt a constant thickness of $0.001\,\text{m}$ in spanwise 
direction. Slip-wall boundary 
conditions are employed for all boundaries except for the inflow, where we prescribe non-reflective 
inflow 
boundary conditions based on Riemann invariants \cite{Computational1992}. Two different fluid mesh resolutions are 
used: 
\textit{Mesh A} contains $123,400$ cells with grid stretching applied in flow direction close to the panel and 
\textit{Mesh B} utilizes a homogeneous grid with $1.82\,$ million cells, see Fig.~\ref{fig:shock_panel_grid}. 
The panel is discretized using $65 \times 2$ ($l = 0.05\,\text{m}$) or $55 \times 2$ ($l = 0.04\,\text{m}$) 
tri-linearly interpolated hexahedral elements. 
For both cases the panel is fully clamped at the bottom, and symmetry boundary conditions are applied in spanwise 
direction. EAS is used in order to avoid shear locking, which may 
affect the solution in such bending-dominated problems when using first-order displacement-based elements. The time
integration factor $\theta = 0.66$ is chosen for the structural time integration. A CFL number of $0.6$ is set for all
simulations.
\begin{figure}
  \centering
 \includegraphics[width= 0.5\columnwidth]{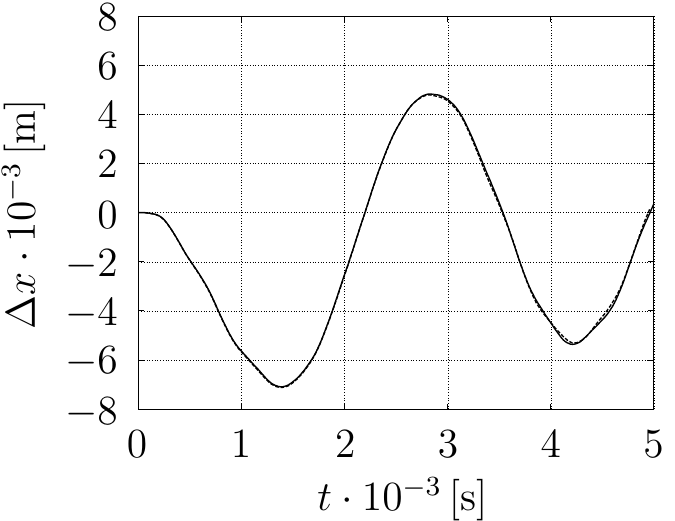}
  \caption{Time evolution of panel tip displacement for $50\,$mm panel length using different fluid meshes. (---) 
\textit{Mesh A}, ($- - -$) \textit{Mesh B}}
  \label{fig:shock_panel_50mm_convergence}  
\end{figure}

Fig.~\ref{fig:shock_panel_50mm_convergence} shows 
the time evolution of the horizontal displacement at the panel tip for the $0.05\,\text{m}$ panel length case on 
\textit{Mesh A} and \textit{Mesh B}. The panel motion is almost identical for both fluid meshes 
throughout the entire simulation time. Results presented below are obtained on fluid 
\textit{Mesh B}.

\begin{figure}
  \centering
 \includegraphics[width= 1.0\columnwidth]{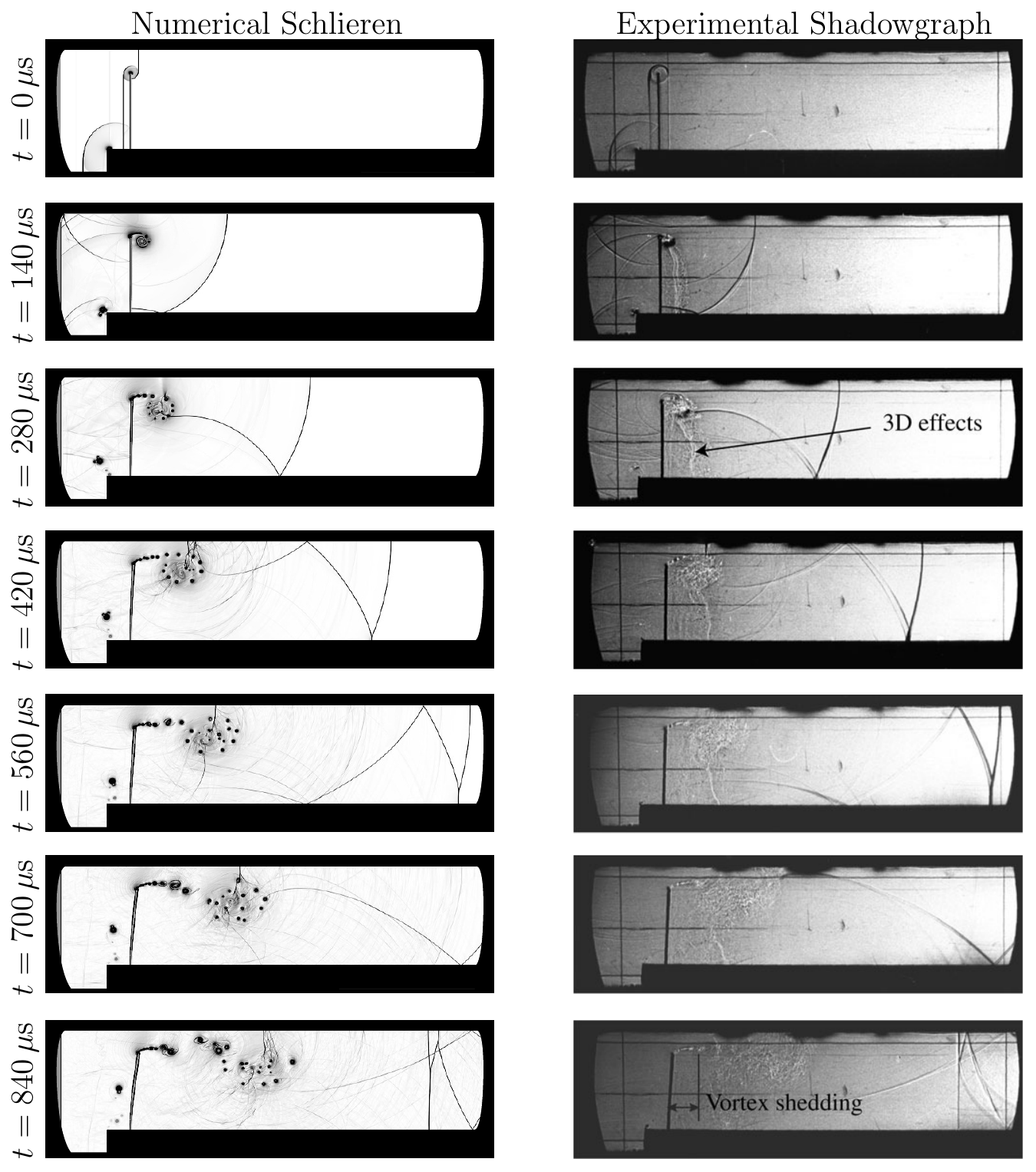}
  \caption{Qualitative comparison between simulation (left) and experiment \cite{Giordano2005} (right) for $50\,$mm 
panel length by means of schlieren images for selected time instances.}
  \label{fig:shock_panel_schlierencomparison}  
\end{figure}

\begin{figure}
  \centering
 \includegraphics[width= 1.0\columnwidth]{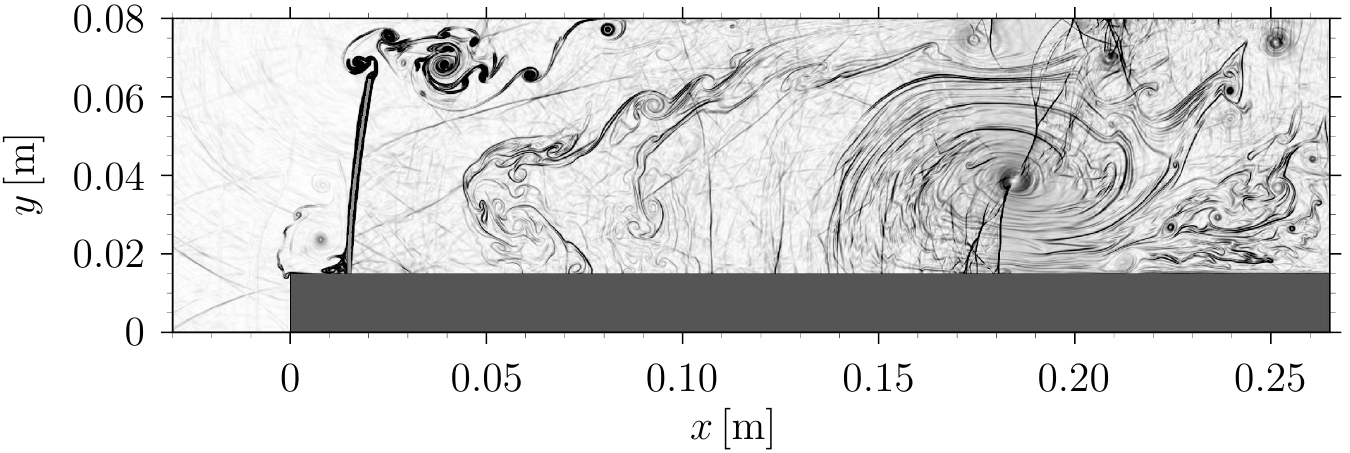}
  \caption{Contour of density gradient magnitude at $t = 4.17\,\text{ms}$.}
  \label{fig:shock_panel_t_add}
\end{figure}

We start with a qualitative analysis of the flow field for the $0.05\,\text{m}$ panel. 
Fig.~\ref{fig:shock_panel_schlierencomparison} shows numerical schlieren (left) and experimental 
shadowgraph visualizations (right) 
extracted from~\cite{Giordano2005} at a time interval of $\Delta t = 
140\,\mu\text{s}$ for a time period of $T = 840\,\mu\text{s}$. At $t = 0\,\mu\text{s}$, the incident right-running 
shock wave has already hit the panel and base plate, leading to the formation of reflected and transmitted shock waves. 
Downstream of the panel the initially normal shock undergoes transition to a cylindrical shock front due to sudden 
area increase ($t = 140\,\mu\text{s}$). While being reflected at the lower wall ($t = 280\,\mu\text{s}$) and traveling 
downstream, it 
undergoes a transition from regular to Mach reflection ($t = 280 - 420\,\mu\text{s}$) and is subsequently reflected at 
the end wall ($t = 700 - 840\,\mu\text{s}$). A main vortex is initially produced at the panel tip due to the 
roll-up of the slipstream accompanied by a vortex shedding process. All flow characteristics described above match 
the experimental results without any notable time lag. However, three-dimensional effects due to leaks between the 
panel and the shock tube side walls are observed in the experiment ($t = 280\,\mu\text{s}$).
Fig.~\ref{fig:shock_panel_t_add} shows a numerical schlieren image at $t = 
4.17\,\text{ms}$, illustrating the maximum panel deflection together with the interaction of the main vortex and the 
upstream moving shock wave.

\begin{figure}
  \centering
 \includegraphics[width= 1.0\columnwidth]{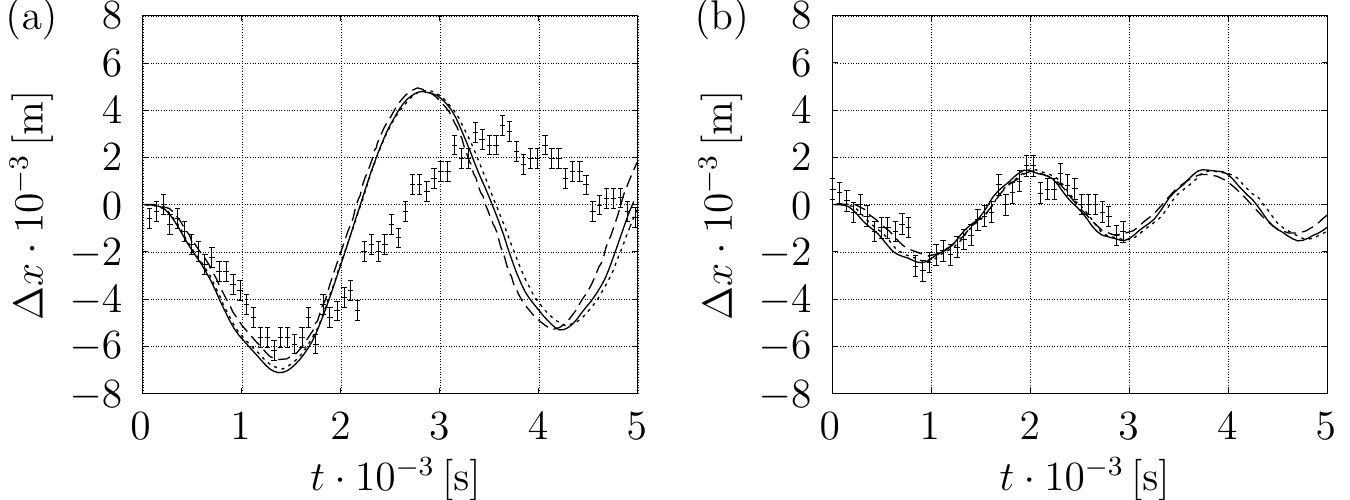}
  \caption{Time evolution of panel tip displacement for (a) $0.05\,$m and (b) $0.04\,$m panel length. (---) present 
results, ($- - -$) \citet{Giordano2005}, ($\cdots$) \citet{Sanches2014}. Error bars denote experimental data 
\cite{Giordano2005}.}
  \label{fig:shock_panel_displ}  
\end{figure}

A quantitative analysis is presented in Fig.~\ref{fig:shock_panel_displ}, where the time evolution of the 
horizontal panel tip displacement is plotted. Fig.~\ref{fig:shock_panel_displ}(a) refers to 
the $0.05\,\text{m}$ panel length case and Fig.~\ref{fig:shock_panel_displ}(b) to the $0.04\,\text{m}$ case, 
respectively. 
In addition to experimental values \cite{Giordano2005} represented through error bars, we include recent inviscid 
numerical results of \citet{Sanches2014}, who employed a finite element based partitioned FSI 
approach utilizing the ALE description to account for moving boundaries and coupling 
with Lagrangian shell elements. Moreover, numerical results by \citet{Giordano2005} are added, who assumed a 
two-dimensional but viscous flow in the laminar regime. For the $0.05\,\text{m}$ panel case, see 
Fig.~\ref{fig:shock_panel_displ}(a), it is observed that all numerical simulations predict a very similar oscillation 
of the panel with respect to the maximum amplitude and frequency of the first period. In comparison to the experimental 
values, both frequency and amplitude of the panel oscillation differ from numerical findings. According to 
\citet{Giordano2005} this difference may be attributed to the lack of damping in the structural model, which, 
however, should be negligible at least for the first period. Another explanation given by the authors relates to small 
deformations of the base in the direct vicinity of the fixing point, which would slightly 
alter both frequency and amplitude of the panel motion. The panel oscillation period obtained with our method is 
$2.85\,\text{ms}$, which is very close to the analytical period of $2.87\,\text{ms}$ when considering the first 
eigenmode of a clamped plate submitted to an impulse load \cite{Giordano2005}. The experimental period is given as 
$3.8\,\text{ms}$.

\begin{figure}
  \centering
 \includegraphics[width= 1.0\columnwidth]{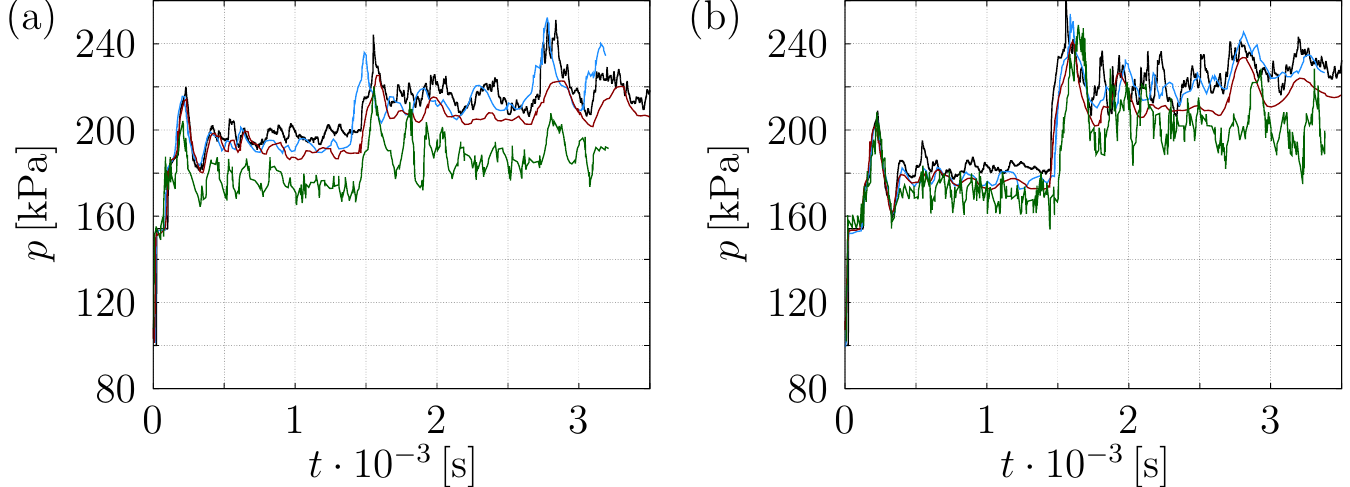}
  \caption{Pressure signal recorded at sensor position (see Fig.~\ref{fig:shock_panel_setup} for exact location of 
pressure probe) for (a) $0.05\,$m and (b) $0.04\,$m panel length. (\textcolor{color1}{\textbf{---}}) present 
results, 
(\textcolor{color2}{\textbf{---}}) \citet{Giordano2005}, (\textcolor{color3}{\textbf{---}}) \citet{Sanches2014}, 
(\textcolor{color4}{\textbf{---}}) experimental values \cite{Giordano2005}. (For interpretation of the references to 
color in this figure legend, the reader is referred to the web version of this article.)}
  \label{fig:shock_panel_pressure} 
\end{figure}

Due to these uncertainties, a second case with $0.04\,\text{m}$ panel 
length has been studied experimentally and numerically in \cite{Giordano2005}. With the shorter panel, the stresses on 
the base part are reduced, which also diminishes the 
influence of the base on the panel motion. We observe excellent agreement with experimental data and numerical 
references, see Fig.~\ref{fig:shock_panel_displ}(b).

Finally, the pressure signals recorded at $\left(x,y \right) = \left(0.035, 0.08\right)\,\text{m}$ 
for both panel lengths are compared to the same numerical and experimental database in 
Fig.~\ref{fig:shock_panel_pressure}. Again, all 
numerical results are similar with respect to the time of arrival of pressure waves at the sensor and the pressure 
difference across the waves. While larger deviations are observed between numerical and experimental data for the 
$50\,\text{mm}$ panel case, almost identical time evolution up to $t = 2\,\text{ms}$ is 
observed for the $0.04\,\text{m}$ panel case. After that time, the pressure obtained experimentally drops 
continuously due to the arrival of reflected expansion waves inside the shock tube, which are not taken into account in 
the numerical simulations.

\subsection{Flutter of a flat plate}
\label{sec:flutter}
\begin{figure}
  \centering
 \includegraphics[width= 0.8\columnwidth]{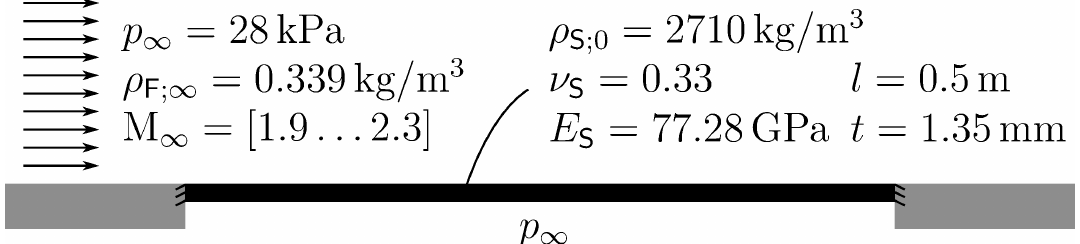}
  \caption{Schematic and main parameters of the flutter problem.}
  \label{fig:flutter_setup}  
\end{figure}

Panel flutter is a self-excited, dynamic aeroelastic instability of thin plate structures, which frequently occurs in 
supersonic flow and is caused by an interaction between aerodynamic, inertial and elastic forces of the system 
\cite{Dowell1972}. For the setup considered here, see Fig.~\ref{fig:flutter_setup}, linear instability theory predicts 
a critical Mach number of $\widetilde{\text{Ma}}_{\idx{crit}} = 2.0$ above which a continuous growth of oscillations 
amplitudes is expected \cite{Dowell1974}. To trigger the instability, the pressure acting on the bottom of the panel 
initially is decreased by $0.1\,\%$ and is kept at this condition for $4\,\text{ms}$. After this time period, the 
pressure is set back to the free-stream pressure. Since the limit Mach number of $\widetilde{\text{Ma}}_{\idx{crit}} = 
2.0$ describes a perfect oscillation without damping or amplification \cite{Piperno1997}, this test case assesses 
effects of numerical damping present in our algorithm.

\begin{figure}
  \centering
 \includegraphics[width= 1.0\columnwidth]{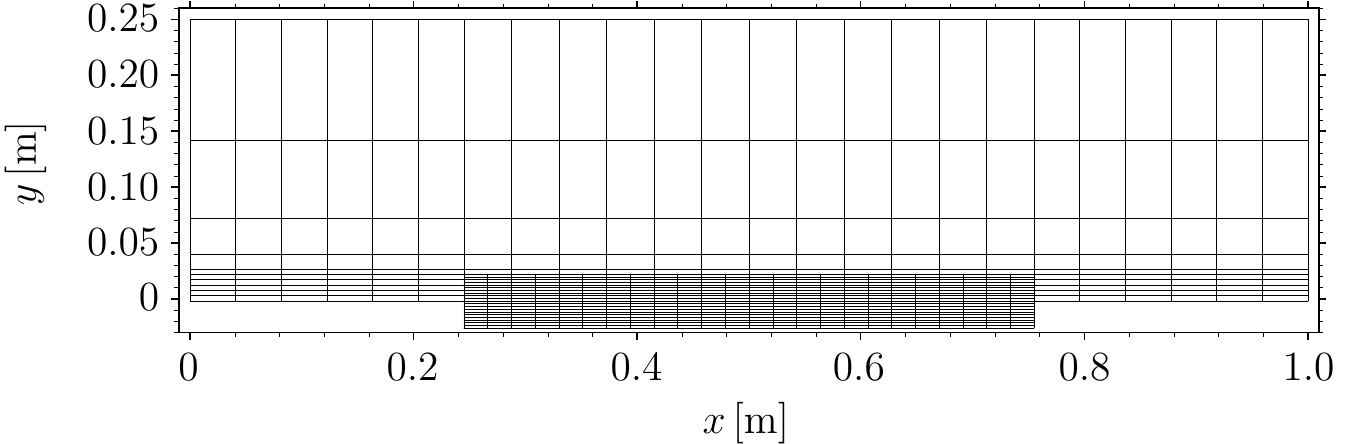}
  \caption{General view of the computational domain and mesh resolution. Every $5$th grid line is shown in the 
$x$- and $y$-direction, respectively.}
  \label{fig:flutter_grid}  
\end{figure}

We consider a supersonic inviscid flow over a flat plate that is clamped at both ends, see
Fig.~\ref{fig:flutter_setup}. The plate has a length of $l = 0.5\,\text{m}$, a thickness of $t = 0.00135\,\text{m}$, a
Young's modulus of $E_{\idx{S}} = 77.28\,\text{GPa}$, a Poisson's ratio of $\nu_{\idx{S}} = 0.33$, and a density of
$\rho_{\idx{S};0} = 2710\,\text{kg/m}^3$. The structure is discretized using $200 \times 8$ tri-linearly interpolated
hexahedral elements in streamwise and wall-normal direction, respectively. To 
avoid shear locking phenomena, the EAS method is used. Results obtained 
with tri-quadratically interpolated hexahedral elements and the same mesh resolution showed only negligible 
differences. If 
not stated otherwise, a geometric linear analysis of the
structure is performed for comparison with references from the literature. The time
integration factor $\theta = 0.5$ is chosen in order to reduce numerical damping. 
The fluid free-stream properties are: $\rho_{\idx{F};\infty} = 0.339\,\text{kg/m}^3$, $p_\infty = 28\,\text{kPa}$ and 
$\text{Ma} = [1.9\ldots2.3]$. 
The computational domain and the fluid mesh
resolution is shown in Fig.~\ref{fig:flutter_grid}. For the results presented here, a grid-converged solution with 
respect to the fluid 
domain has been obtained with a total number of
$16,500$ cells. The grid is uniform in the region around the panel ($0.25\,\text{m}
\le x \le 0.75\,\text{m}$) with a cell size of $\Delta x = 4.25\times 10^{-3}\,\text{m}$ and $\Delta y = 4.8\times
10^{-4}\,\text{m}$. A cavity of height $h = 2.2\times 10^{-2}\,\text{m}$ is added below the panel ($y \le 0\,\text{m}$) 
to
account for the panel motion in this region. Since the problem is two-dimensional, we adopted a 
constant thickness of
$\Delta z = 5\times 10^{-3}\,\text{m}$ in spanwise direction. Slip-wall boundary conditions are imposed at all boundary 
patches except for
the inflow and outflow patch. At the inflow we prescribe all flow quantities which leads to a fully reflective boundary
condition. At the outflow we perform linear extrapolation. The CFL number is $0.6$ for all simulations.

\begin{figure}
  \centering
  \includegraphics[width= 0.8\columnwidth]{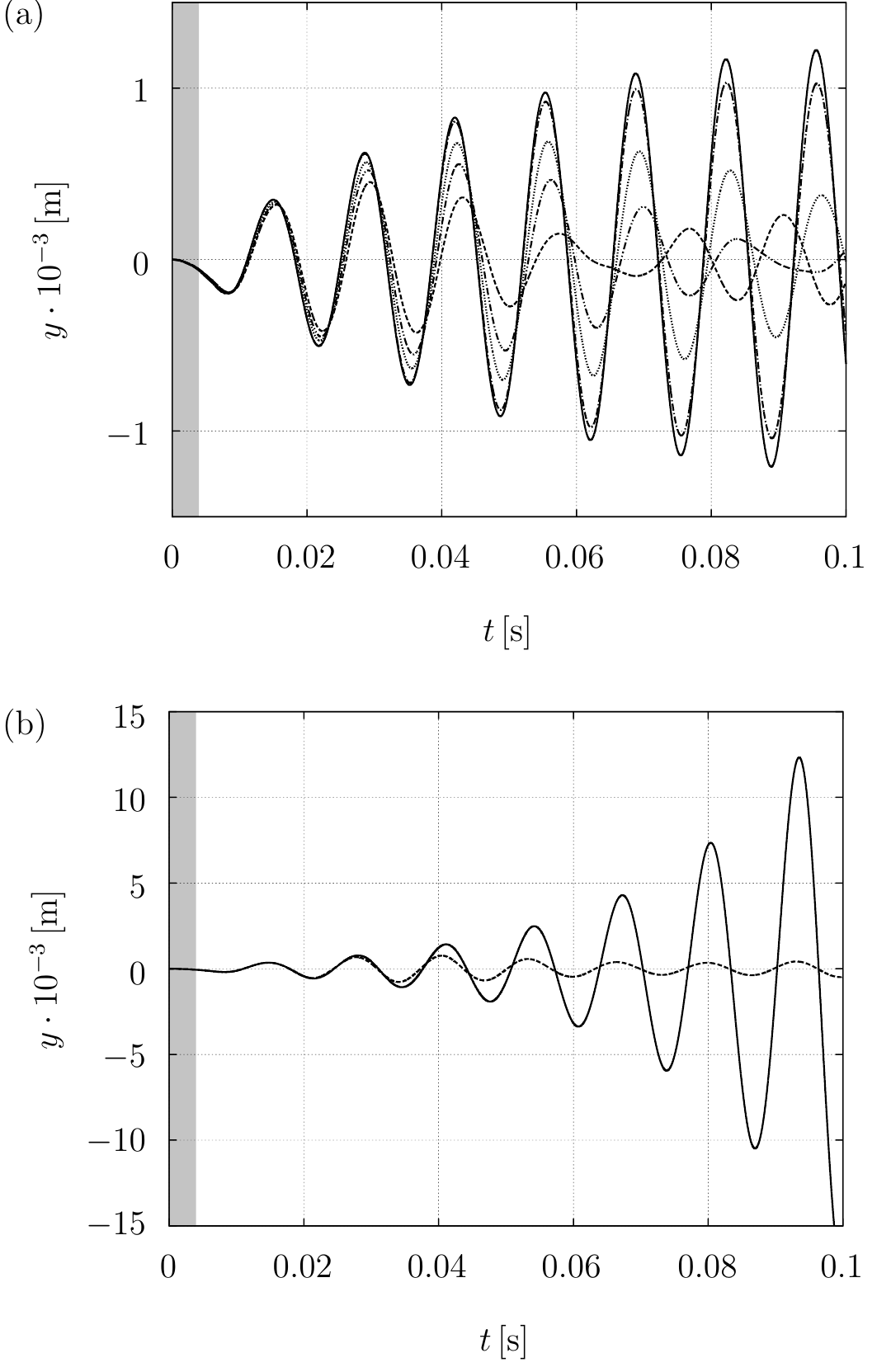}
  \caption{(a) Vertical deflection of the plate at $x = 0.6\,\text{m}$ for $\text{Ma} \in [1.9,2.0,2.05,2.09,2.1]$. (- 
- -) $\text{Ma} = 1.9$, ($- \cdot \cdot\, -$) $\text{Ma} = 2.0$, (\,$\cdots$) $\text{Ma} = 2.05$, ($- \cdot -$) 
$\text{Ma} = 2.09$, (---) $\text{Ma} = 2.1$. (b) Geometrically linear and nonlinear plate deflections at $x = 
0.6\,\text{m}$ for $\text{Ma} = 2.3$\,. 
(---) linear, (- - -) nonlinear. The gray shaded area indicates the initial perturbation time.}
  \label{fig:flutter_displ_M_all} 
\end{figure}

The time evolution of the vertical displacement of the panel at the streamwise position $x = 0.6\,\text{m}$ for Mach 
numbers $\text{Ma} = \left[1.9, 2.0, 2.05, 2.09, 2.1\right]$ is shown in Fig.~\ref{fig:flutter_displ_M_all}(a). The gray 
shaded area indicates the initial perturbation time. 
While the panel oscillations for Mach numbers below $\text{Ma} = 2.09$ are damped, amplification of 
panel deflection can be observed for $\text{Ma} = 2.1$. We found the limit Mach number to be 
$\text{Ma}_{\idx{crit}} = 2.09$, which is close to the analytical solution 
($\widetilde{\text{Ma}}_{\idx{crit}} = 2.0$) with an error of $4.5\%$ and to numerical results reported by 
\citet{Teixeira2005} and \citet{Sanches2014} ($\text{Ma}_{\idx{crit}} = 2.05$). Fig.~\ref{fig:flutter_displ_M_all}(b) 
shows a comparison between geometrically linear and nonlinear panel solutions for a Mach number of $\text{Ma} = 2.3$. 
Exponential growth of the initial disturbance is observed for linear theory, which confirms analytical and 
numerical results \cite{Dowell1974, Piperno1997, Teixeira2005, Sanches2014}. In the geometrically nonlinear 
case, limited displacement amplitudes are observed. According to \citet{Dowell1972}, the behavior 
of the panel after flutter onset is mainly dominated by structural nonlinearities. Nonlinear structural coupling 
between bending and stretching of the plate may in fact increase its effective stiffness, thereby modifying the 
dynamic response of the system.

\begin{figure}
  \centering
  \includegraphics[width= 1.0\columnwidth]{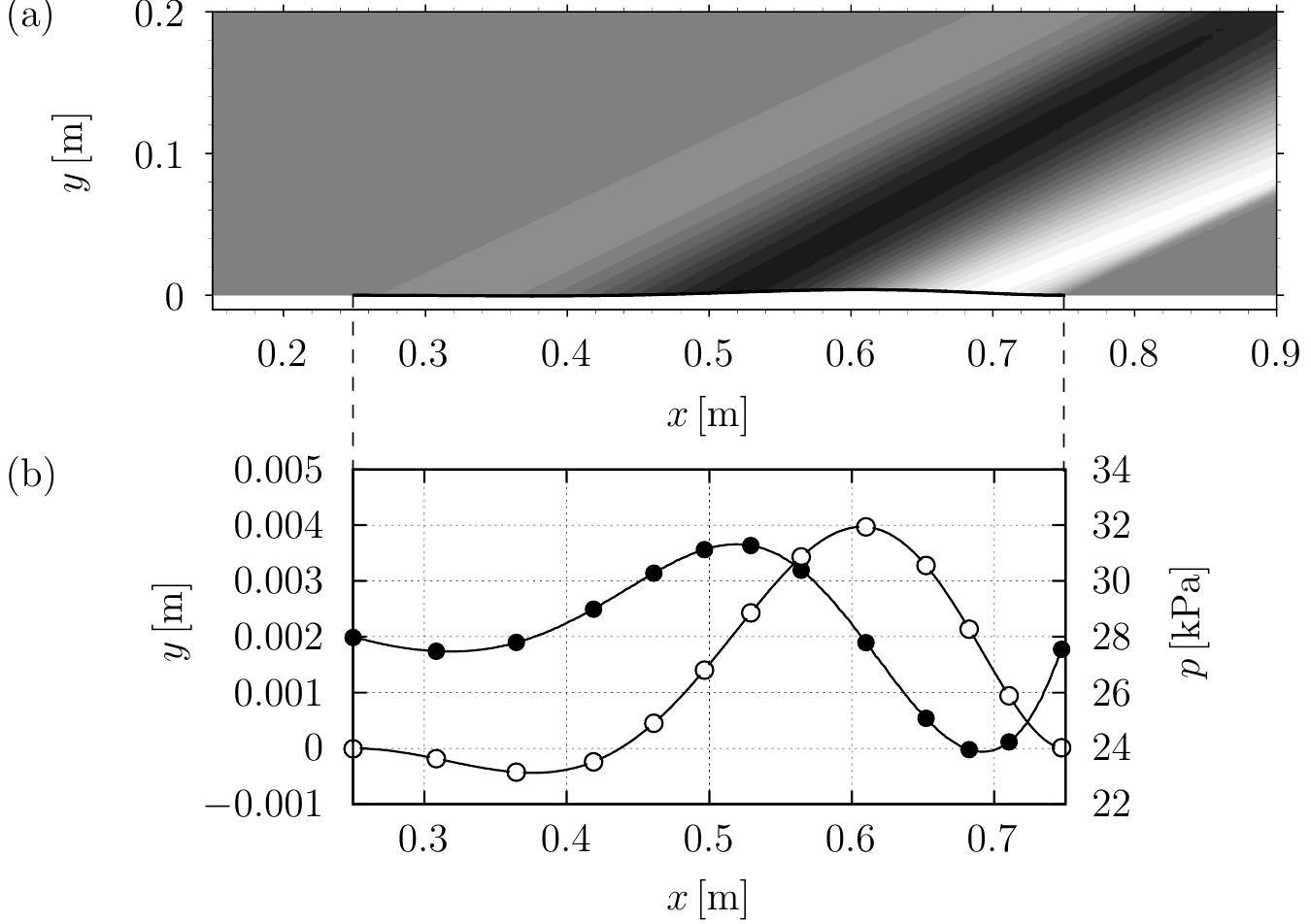}
  \caption{(a) Pressure distribution for $\text{Ma} = 2.3$ at $t = 0.068\,\text{s}$. Color scale from white to black 
using $20$ equally spaced contour levels for $p \in [24-32]\,\text{kPa}$. (b) Interface 
pressure and associated 
flutter mode. \mbox{(--- $\bullet$ ---) $p^{\Gamma}$}, \mbox{(--- $\circ$ ---) $\Delta y$}.  }
  \label{fig:flutter_M_2_3_linear_pressure} 
\end{figure}

Finally, the pressure distribution in the fluid domain together with the associated flutter mode and 
wall-pressure distribution at time instant $t = 0.068\,\text{s}$ is shown in 
Fig.~\ref{fig:flutter_M_2_3_linear_pressure}. The deflection of the panel leads to the formation of compression and 
expansion waves in the fluid. Compression waves are observed for a positive interface slope, whereas expansion waves 
occur for negative interface slopes, which is consistent with Ackeret's linear theory. The maximum displacement for the 
flutter mode is found at $70\%$ of the panel length, confirming analytical \cite{Houbolt1958, Dowell1974} and numerical 
\cite{Piperno1997, Teixeira2005, Sanches2014} findings. Local minima and maxima in the 
wall-pressure distribution in Fig.~\ref{fig:flutter_M_2_3_linear_pressure}(b) coincide with interface 
inflection points. The smooth wall-pressure distribution confirms once again the accurate interface treatment.

\subsection{Grid convergence study}
\label{sec:accuracy}

The accuracy of the computed solution is verified through a grid convergence study. The simulation setup is similar to 
the case presented in Section~\ref{sec:shock_cylinder}.
The formerly rigid cylinder is now replaced by an elastic structure and the wind-tunnel walls are removed.  
For the cylinder, which is initially located at $(x,y) = (0.15,0.0)\,\text{m}$, a Young's modulus of $E_{\idx{S}} = 
800\,\text{Pa}$, a Poisson's ratio of $\nu_{\idx{S}} = 0.3$, and a density of $\rho_{\idx{S};0} = 15\,\text{kg/m}^3$ 
have been adopted. The remaining parameters are identical to the setup described in Section~\ref{sec:shock_cylinder}. 

\begin{figure}[htbp]
\centering
\includegraphics[width= 1.0\columnwidth]{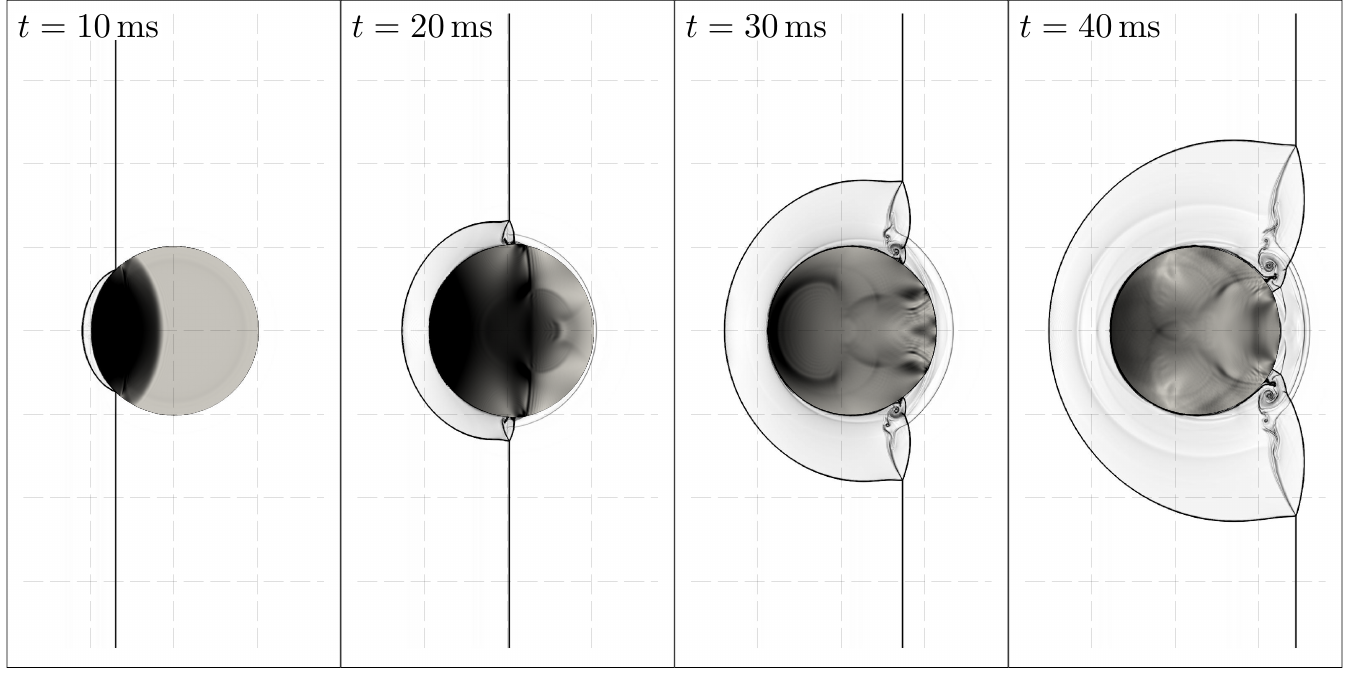}
  \caption{Contours of density gradient magnitude in the fluid domain and magnitude of the Cauchy stress tensor in the
solid domain at four different simulation times for the reference grid $\mathcal{G}^{ref}$. }
  \label{fig:shock_deformable_cylinder_grid_ref}
\end{figure}

Since no analytical solution for this complex interaction exists, we have performed a well resolved reference 
simulation. 
The reference grid, in the following denoted as $\mathcal{G}^{ref}$, has a resolution of $1280 \times 2560$ 
cells in the fluid domain and spatial dimensions of $0.2\,\text{m} \times 0.4\,\text{m}$. The 
cylinder is discretized with $2048$ tri-linearly interpolated hexahedral elements along its circumference. For the 
remaining grids $\left .\mathcal{G}^{k}\right|_{k = 1 \ldots 5}$, where $\mathcal{G}^{5}$ denotes the finest grid, the
fluid resolution is successively halved and the unstructured mesh resolution of the solid is halved in radial and
circumferential direction. A uniform time step of $\Delta t = 5.1\times 
10^{-6}\,\mathrm{s}$ is used for all simulations, which corresponds to the maximum allowable time step size for the 
reference simulation at a CFL number of $0.6$.

\begin{figure}[htbp]
\centering 
\includegraphics[width=0.7\columnwidth]{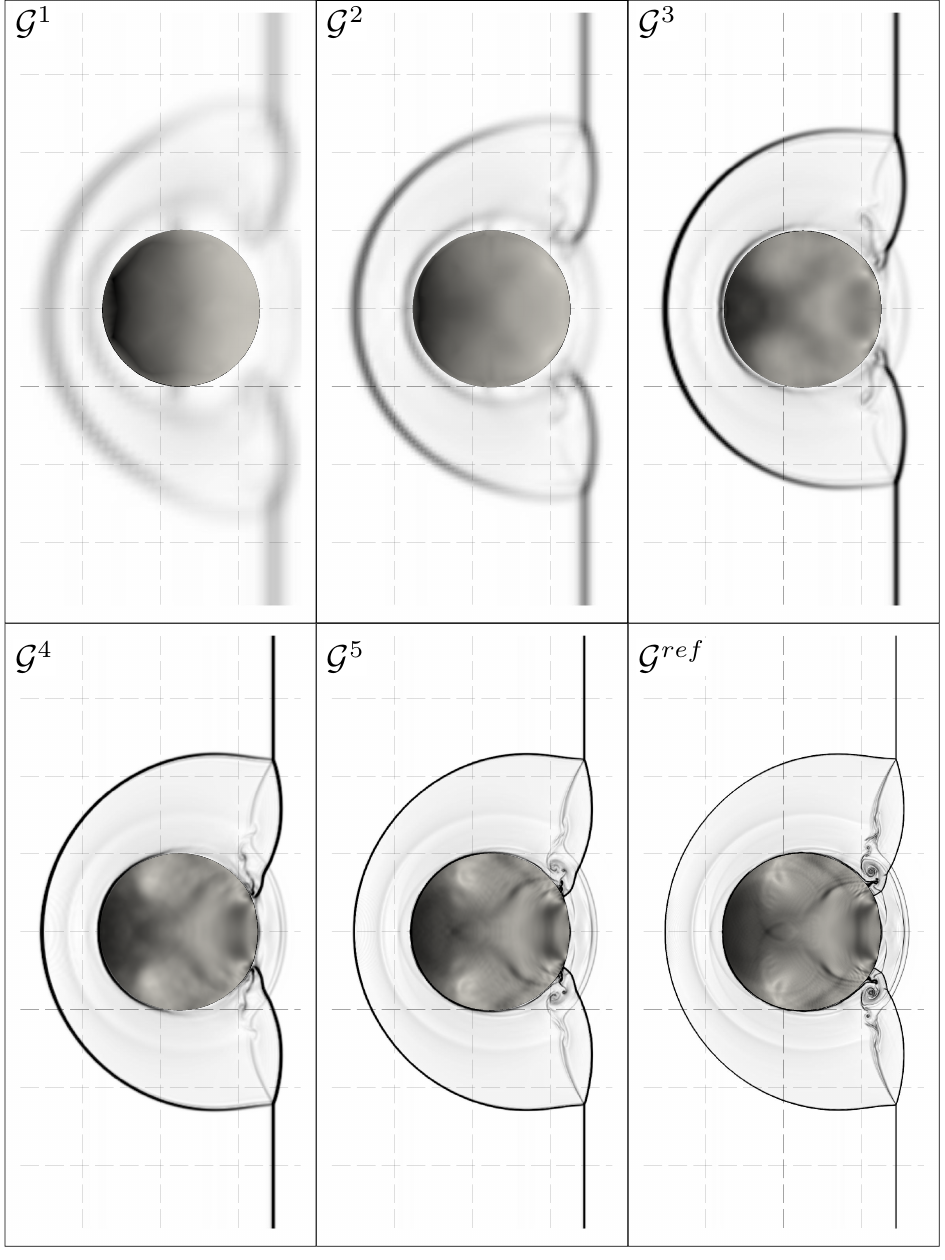}
  \caption{Contours of density gradient magnitude in the fluid and magnitude of the Cauchy stress tensor in the
solid domain at time $t = 40\,\mathrm{ms}$ for all considered mesh resolutions.}
  \label{fig:shock_deformable_cylinder_all_grids}
\end{figure}

Figure~\ref{fig:shock_deformable_cylinder_grid_ref} shows a numerical schlieren visualization of the resulting flow 
field together with 
the Cauchy stress field within the solid at times $t = 10,\,20,\,30$ and $40\,\mathrm{ms}$ computed on 
$\mathcal{G}^{ref}$. As expected, both fields are symmetric with respect to the $x$ axis, even though no symmetry is 
presumed for the algorithm. At time $t = 10\,\mathrm{ms}$, the 
incident shock has already hit the cylinder and is subsequently reflected. The impact on the cylinder generates a
shock wave which propagates through the solid. As the shock travels further around the 
cylinder, it undergoes transition from regular to Mach reflection ($t = 20\,\mathrm{ms}$). At the same time, the 
windward side of the cylinder is compressed, while the leeward side moves slightly downstream and generates a shock 
wave in the fluid. At the triple point, which connects the incident shock, the reflected shock and the Mach stem, a 
contact discontinuity develops. By the time the cylinder is accelerated ($t = 30\,\mathrm{ms}$), the reflected shock 
has propagated further upstream and a roll-up of the contact discontinuity is observed, which is enhanced by the 
interaction with the leeward shock wave. At the final time $t = 40\,\mathrm{ms}$, several shock waves emerging from the 
fluid-structure interface can be observed and an overall complex flow field has developed.

As a qualitative measure of the accuracy of our method, Fig.~\ref{fig:shock_deformable_cylinder_all_grids} shows 
numerical schlieren visualizations together with Cauchy stresses at the final time $t = 40\,\mathrm{ms}$ for all mesh 
resolutions $\mathcal{G}^{k}$. While the 
overall results with respect to the final cylinder position and the incident and reflected shock wave within the 
fluid domain agree well between all mesh resolutions, the finer grids ($\mathcal{G}^{3}, \mathcal{G}^{4}, 
\mathcal{G}^{5}$) provide fine scale features such as contact discontinuities and weak shock waves emerging from the 
cylinder surface which are partially missing or not well resolved on the coarse grids  ($\mathcal{G}^{1}, 
\mathcal{G}^{2}$).

\ctable[
        cap     = {Computed $L_{1}$ and $L_{2}$ fluid error norms $\mathcal{E}_{1}^{k}$ and 
$\mathcal{E}_{2}^{k}$ on all grids $\mathcal{G}^{k}$ with respect to density, pressure and velocity magnitude at time 
$t 
= 
40\,\mathrm{ms}$. Estimated convergence rates $m$ are based on a least squares fit.},
        caption = {Computed $L_{1}$ and $L_{2}$ fluid error norms $\mathcal{E}_{1}^{k}$ and 
$\mathcal{E}_{2}^{k}$ on all grids $\mathcal{G}^{k}$ with respect to density, pressure and velocity magnitude at time 
$t = 40\,\mathrm{ms}$. Estimated convergence rates $m$ are based on a least squares fit.},
        label= {table:spatial_accuracy_fluid},
        pos = {t},
        maxwidth = \columnwidth,
        captionskip = -0.4cm,
        doinside = \tiny,
        notespar,
        nosuper,
 ]{lrrrrrrrr}{
           \tnote[$a.$]{Fluid cell size in $[\mathrm{m}]$. A uniform grid is used.}
           \tnote[$b.$]{Structural element length along the cylinder circumference in $[\mathrm{m}]$.}
        %
}{
\FL
$\mathcal{G}^{k}$&
$\Delta_{\idx{F}}$\tmark[$a$] &
$\Delta_{\idx{S}}$\tmark[$b$] &
$\mathcal{E}_{1,\rho}^{k}$& 
$\mathcal{E}_{2,\rho}^{k}$& 
$\mathcal{E}_{1,p}^{k}$& 
$\mathcal{E}_{2,p}^{k}$& 
$\mathcal{E}_{1,\left|\discret{u}\right|}^{k}$& 
$\mathcal{E}_{2,\left|\discret{u}\right|}^{k}$ \\
\FL
$\mathcal{G}^{1}$ & $5\cdot 10^{-3}$     & $4.91\cdot 10^{-3}$    & $2.1\cdot 10^{-1}$   & $4.3\cdot 10^{-1}$
                                                                  & $7.8\cdot 10^{-1}$   & $1.8\cdot 10^{ 0}$  
                                                                  & $1.2\cdot 10^{-1}$   & $3.0\cdot 10^{-1}$\\
$\mathcal{G}^{2}$ & $2.5\cdot 10^{-3}$   & $2.45\cdot 10^{-3}$    & $8.5\cdot 10^{-2}$   & $2.2\cdot 10^{-1}$   
                                                                  & $3.2\cdot 10^{-1}$   & $8.8\cdot 10^{-1}$ 
                                                                  & $4.6\cdot 10^{-2}$   & $1.7\cdot 10^{-1}$\\
$\mathcal{G}^{3}$ & $1.25\cdot 10^{-3}$  & $1.23\cdot 10^{-3}$    & $2.7\cdot 10^{-2}$   & $8.9\cdot 10^{-2}$  
                                                                  & $9.7\cdot 10^{-2}$   & $3.3\cdot 10^{-1}$         
                                                                  & $1.9\cdot 10^{-2}$   & $9.5\cdot 10^{-2}$\\
$\mathcal{G}^{4}$ & $6.25\cdot 10^{-4}$  & $6.14\cdot 10^{-4}$    & $1.3\cdot 10^{-2}$   & $5.9\cdot 10^{-2}$          
                                                                  & $4.3\cdot 10^{-2}$   & $2.1\cdot 10^{-1}$          
                                                                  & $8.4\cdot 10^{-3}$   & $6.0\cdot 10^{-2}$\\
$\mathcal{G}^{5}$ & $3.125\cdot 10^{-4}$ & $3.07\cdot 10^{-4}$    & $5.2\cdot 10^{-3}$   & $3.0\cdot 10^{-2}$          
                                                                  & $1.6\cdot 10^{-2}$   & $1.0\cdot 10^{-1}$           
                                                                  & $3.5\cdot 10^{-3}$   & $3.0\cdot 10^{-2}$\\
\cmidrule(l){4-9}
& & \textsc{Rate $m$} & $1.36$ & $1.03$ & $1.37$ & $1.12$ & $1.34$ & $0.80$
\LL}

\ctable[
        cap     = {Computed $L_{1}$ and $L_{2}$ structural error norms $\mathcal{E}_{1}^{k}$ and 
$\mathcal{E}_{2}^{k}$ with respect to the interface displacement magnitude and $\mathcal{E}^{k}$ with
respect to the interface force in $x$-direction on all grids $\mathcal{G}^{k}$ at time $t =
40\,\mathrm{ms}$. Estimated convergence rates $m$ are based on a least squares fit.},
        caption = {Computed $L_{1}$ and $L_{2}$ structural error norms $\mathcal{E}_{1}^{k}$ and 
$\mathcal{E}_{2}^{k}$ with respect to the interface displacement magnitude and $\mathcal{E}^{k}$ with
respect to the interface force in $x$-direction on all grids $\mathcal{G}^{k}$ at time $t =
40\,\mathrm{ms}$. Estimated convergence rates $m$ are based on a least squares fit.},
        label= {table:spatial_accuracy_struct},
        pos = {t},
        maxwidth = \columnwidth,
        captionskip = -0.4cm,
        doinside = \small,
        notespar,
        nosuper,
 ]{lrrrrr}{
           \tnote[$a.$]{Fluid cell size in $[\mathrm{m}]$. A uniform grid is used.}
           \tnote[$b.$]{Structural element length along the cylinder circumference in $[\mathrm{m}]$.}
        %
}{
\FL
$\mathcal{G}^{k}$&
$\Delta_{\idx{F}}$\tmark[$a$] &
$\Delta_{\idx{S}}$\tmark[$b$] &
$\mathcal{E}_{1,\left|\discret{d}^\Gamma\right|}^{k}$& 
$\mathcal{E}_{2,\left|\discret{d}^\Gamma\right|}^{k}$&
$\mathcal{E}_{\left| \discret{\sigma}^\Gamma_\idx{F} \cdot \discret{n}^\Gamma_\idx{F} \right|}^{k}$ \\
\FL
$\mathcal{G}^{1}$ & $5\cdot 10^{-3}$     & $4.91\cdot 10^{-3}$    & $8.9\cdot 10^{-4}$   & $7.7\cdot 10^{-4}$
                                                                  & $5.6\cdot 10^{-3}$\\
$\mathcal{G}^{2}$ & $2.5\cdot 10^{-3}$   & $2.45\cdot 10^{-3}$    & $6.1\cdot 10^{-4}$   & $5.4\cdot 10^{-4}$   
                                                                  & $2.0\cdot 10^{-3}$\\
$\mathcal{G}^{3}$ & $1.25\cdot 10^{-3}$  & $1.23\cdot 10^{-3}$    & $1.5\cdot 10^{-4}$   & $1.3\cdot 10^{-4}$  
                                                                  & $5.8\cdot 10^{-4}$\\
$\mathcal{G}^{4}$ & $6.25\cdot 10^{-4}$  & $6.14\cdot 10^{-4}$    & $2.6\cdot 10^{-5}$   & $2.4\cdot 10^{-5}$          
                                                                  & $2.2\cdot 10^{-4}$\\
$\mathcal{G}^{5}$ & $3.125\cdot 10^{-4}$ & $3.07\cdot 10^{-4}$    & $1.8\cdot 10^{-5}$   & $1.6\cdot 10^{-5}$          
                                                                  & $4.0\cdot 10^{-5}$\\
 \cmidrule(l){4-6}
& & \textsc{Rate $m$} & $1.58$ & $1.56$ & $1.75$
\LL}

\begin{figure}[htbp]
  \centering
  \begin{subfigure}{\textwidth}
	  \centering
	  \includegraphics{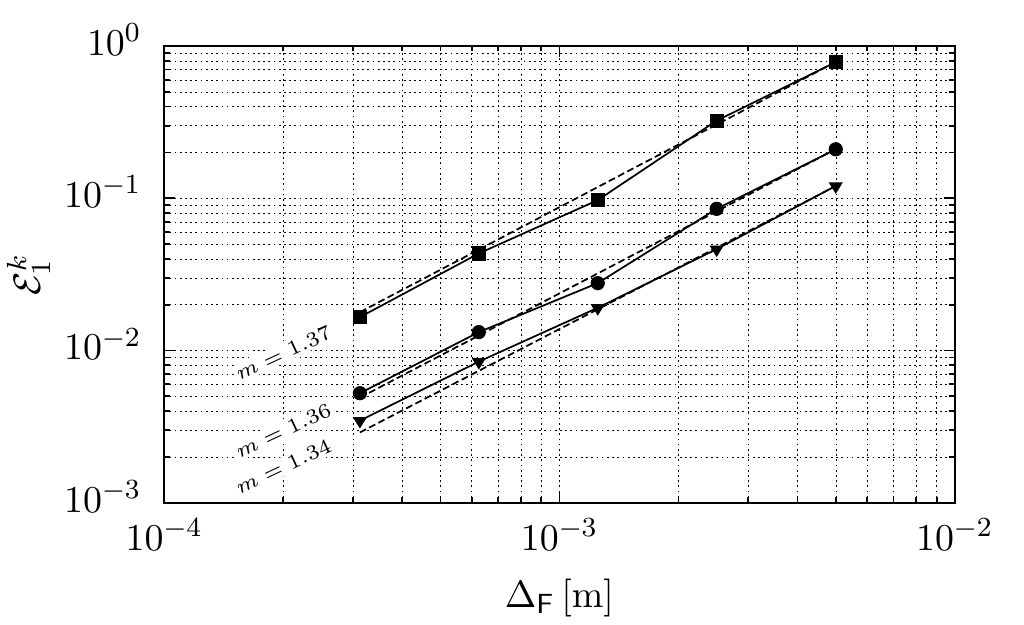}
  \end{subfigure}%
  \\
  \begin{subfigure}{\textwidth}
	  \centering
	  \includegraphics{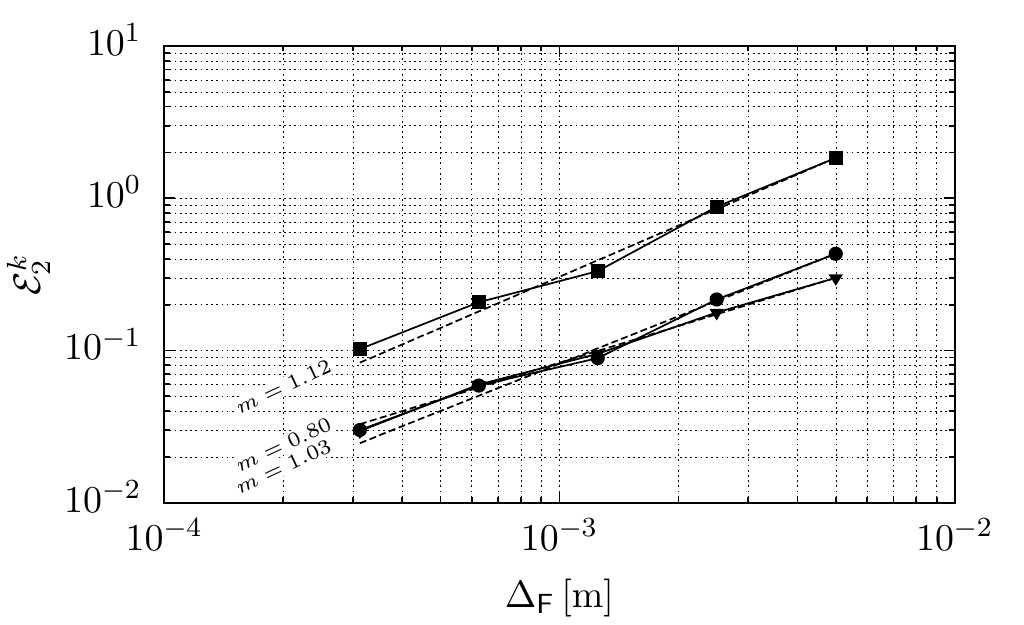}
  \end{subfigure}
  \caption{Computed $L_{1}$ (top) and $L_{2}$ (bottom) fluid error norms $\mathcal{E}_{1}^{k}$ and 
$\mathcal{E}_{2}^{k}$ on all grids $\mathcal{G}^{k}$ at time $t = 40\,\mathrm{ms}$. (--- $\bullet$ ---) 
$\mathcal{E}_{\rho}^{k}$, (--- \solidsquare ---) $\mathcal{E}_{p}^{k}$, (--- $\blacktriangledown$ ---) 
$\mathcal{E}_{\left|\discret{u}\right|}^{k}$. Dashed lines represent least squares fits. 
Estimated convergence rates $m$ are highlighted.}
  \label{fig:spatial_accuracy_fluid}
\end{figure}

\begin{figure}[htbp]
  \centering
    \centering
    \includegraphics{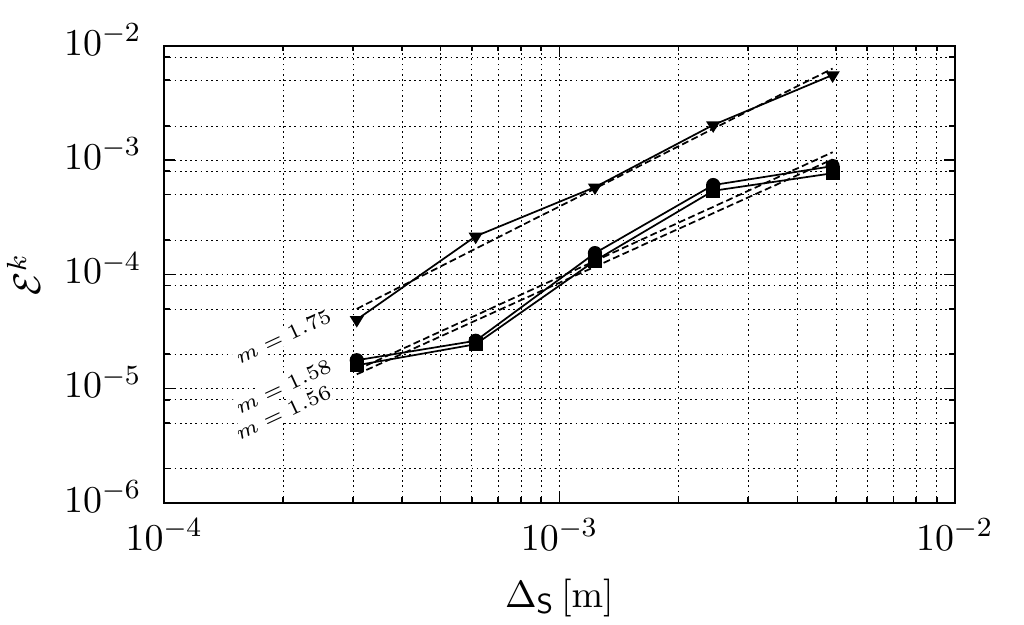}
  \caption{Computed interface norms on all grids $\mathcal{G}^{k}$ at time $t =
40\,\mathrm{ms}$. (--- $\bullet$ ---) $\mathcal{E}_{1,\left|\discret{d}^\Gamma\right|}^{k}$, (--- \solidsquare ---)
$\mathcal{E}_{2,\left|\discret{d}^\Gamma\right|}^{k}$, (--- $\blacktriangledown$ ---)
$\mathcal{E}_{\left| \discret{\sigma}^\Gamma_\idx{F} \cdot \discret{n}^\Gamma_\idx{F} \right|}^{k}$. Dashed lines
represent least squares fits. 
Estimated convergence rates $m$ are highlighted.}
  \label{fig:spatial_accuracy_struct}
\end{figure}

A quantitative measure of accuracy within the fluid domain is given by the 
discrete $L_{p}$ norm of the error for a solution variable $\mathcal{S}^{k}$ on grid $\mathcal{G}^{k}$, which we define 
as
\begin{equation}
 \mathcal{E}_{p}^{k} = \left[ \dfrac{1}{N} \sum_{i=1}^{N} \left( \mathcal{S}_{i}^{k} - \mathcal{S}_{i}^{ref} 
\right )^{p} \right]^{\frac{1}{p}}\, .
 \label{eq:lp-norm}
\end{equation}
Here, $N$ denotes the total number of fluid cells considered on grid $\mathcal{G}^{k}$. In order to evaluate the
convergence of the coupling problem, interface quantities are used. The error of the magnitude of the interface
displacement of the structure is measured by integrating the error over the coupling surface. Thus, the sum in
\eqref{eq:lp-norm} is replaced by an integration and division by the number of grid points is replaced by division by
the area of the coupling surface. A second interface quantity of interest is the coupling force in $x$-direction. The
corresponding error is computed as 
\begin{equation}
 \mathcal{E}^{k} = \left| \int_{\Gamma_{\idx{F}}^k} \discret{\sigma}^{\Gamma,k}_\idx{F} \cdot
\discret{n}^{\Gamma,k}_\idx{F} \mathrm{d} \Gamma - \int_{\Gamma_{\idx{F}}^{ref}} \discret{\sigma}^{\Gamma,ref}_\idx{F}
\cdot \discret{n}^{\Gamma,ref}_\idx{F} \mathrm{d} \Gamma \right| \, .
 \label{eq:force-norm}
\end{equation}
Table~\ref{table:spatial_accuracy_fluid} summarizes estimated errors in the fluid density, pressure and velocity 
magnitude. Table~\ref{table:spatial_accuracy_struct} contains estimated errors in the structural interface 
displacement magnitude and in the coupling force in $x$-direction. Both tables include associated
convergence 
rates at time $t = 40\,\mathrm{ms}$ which are estimated from a 
least squares fit to the logarithm of the errors with the target function $\mathcal{F} = C \cdot
\Delta^{m}$, 
where $\Delta$ denotes either the discrete fluid or structural mesh resolution and $C$ denotes a positive constant
independent of 
the grid. We observe convergence rates with respect to the $L_{1}$ norm of approximately $1.3$ for all fluid variables, 
while the $L_{2}$ convergence rates are overall lower. Similar results have been observed by \citet{Henshaw2008} for a 
pure fluid simulation of shock diffraction by a sphere. As expected from the flow field at time $t = 40\,\mathrm{ms}$, 
which is dominated by shock waves and contact discontinuities, the convergence order with respect to all fluid 
variables is first order. Figure~\ref{fig:spatial_accuracy_fluid} shows the associated variation of the $L_{1}$ and 
$L_{2}$ error norms in the fluid variables plotted against the fluid mesh resolution $\Delta x = \Delta y = 
\Delta_{\idx{F}}$. The error of the magnitude of the interface displacement is of order $1.58$ and
$1.56$ in $L_{1}$ and $L_{2}$ norm, respectively. The convergence rate with respect to coupling force in
$x$-direction computed at the interface from the fluid is $1.75$. Due to the fluid-structure coupling and the use of 
tri-linear finite elements, the overall expected order of convergence is at most second order. 
Figure~\ref{fig:spatial_accuracy_struct} contains the error norms for the magnitude of the interface displacements as 
well as for the coupling force in $x$-direction plotted against the structural mesh resolution $\Delta_{\idx{S}}$ in 
circumferential direction.

This convergence study still has its limitations. First of all the almost standard limitation in
such cases is not to appropriately take into account the coupling of spatial and temporal error but then comparing
spatial errors at a certain point in time. An additional limitation in this case is that the specific FSI example does
not include real structural dynamics in terms of large deformations, but rather shows a combination of rigid body 
dynamics combined with wave
propagation in the solid, which has obviously different features. Given the lack of an established benchmark
example we intended to stay close to a widely accepted example, namely the shock wave impact on a rigid cylinder as
given before. Besides all the limitations, however, the provided convergence study should give some insight into the 
performance of the coupling approach.

\newpage
\section{Numerical Example - Buckling of a three-dimensional inflated thin shell}
\label{sec:numerical_example}

\begin{figure}
  \centering
  \includegraphics[width= 1.\columnwidth]{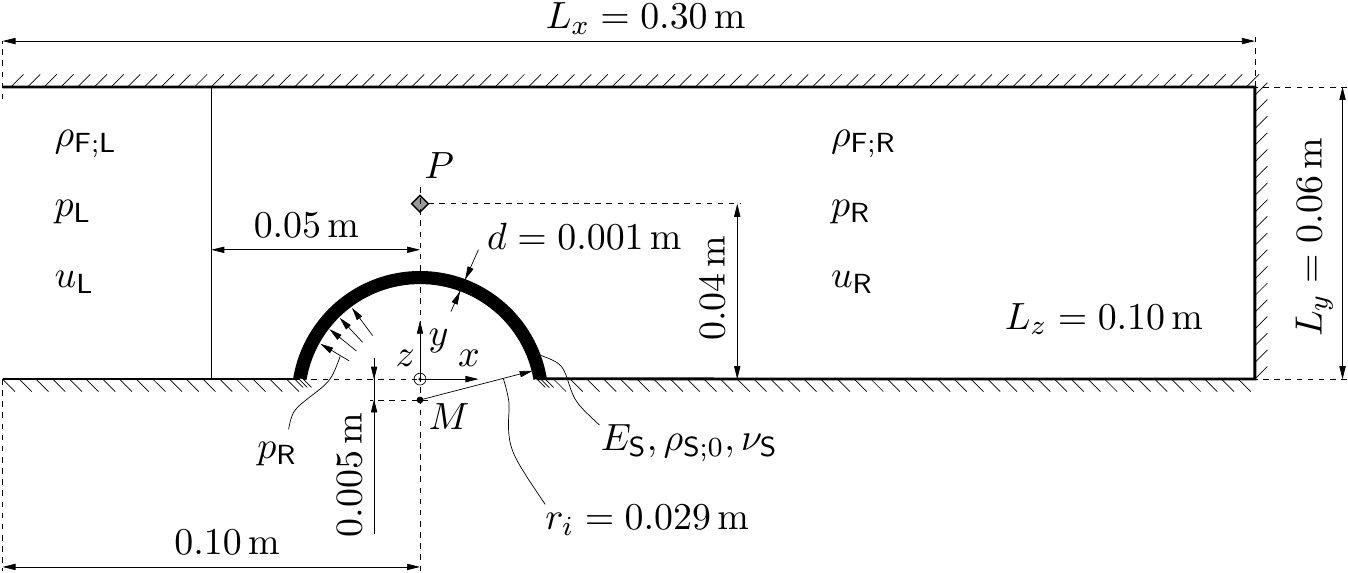}
  \caption{Setup at $xy$-midplane for shock wave impact on a thin-walled shell including geometric dimensions.}
  \label{fig:shock_membrane_setup}
\end{figure}
\begin{figure}
  \centering
  \includegraphics[width= 1.\columnwidth]{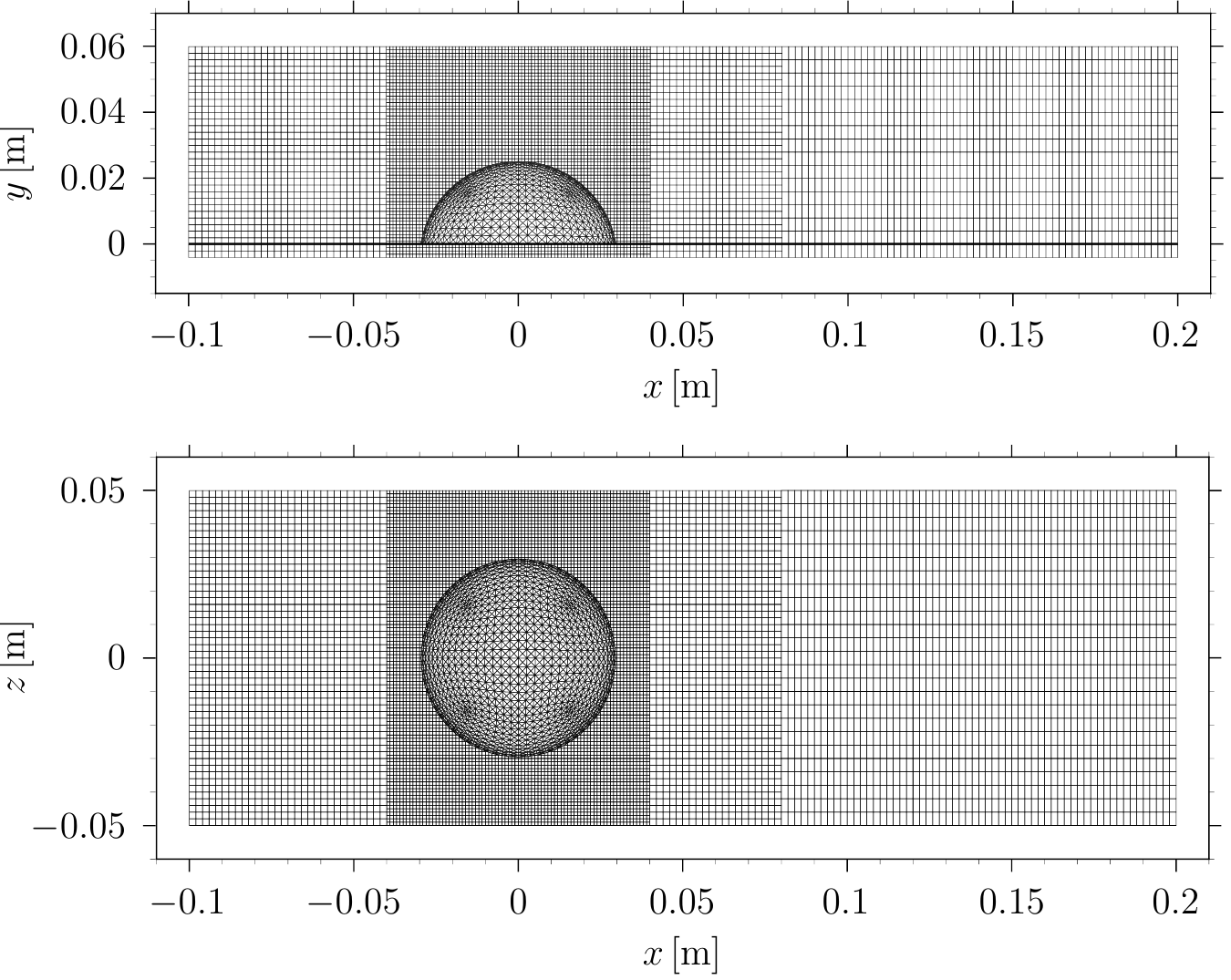}
  \caption{General view of the computational domain and mesh resolution. The triangulated solid interface is  
additionally illustrated.}
  \label{fig:shock_membrane_grid}
\end{figure}
We present a numerical example to show the ability of our method to handle large and complex 
structural deformations in FSI problems. The presented example studies the interaction between a flexible inflated 
thin shell and a $\text{Ma} = 1.21$
shock wave. Pre- and post-shock fluid states are equal to the conditions introduced in  
Section \ref{sec:shock_panel}, with the initial shock position being located at $x = -0.05\,\text{m}$. Details of the 
setup are shown in Fig.~\ref{fig:shock_membrane_setup}. The spherical membrane has a thickness of $d = 
0.001\,\text{m}$ and an inner radius $r_i = 0.029\,\text{m}$ with its
center $M$ located at $(x,y,z) = (0,-0.005,0)\,\text{m}$.
Material properties are $E_{\idx{S}} = 0.07\,\text{GPa}$, $\rho_{\idx{S};0} = 1000\,\text{kg/m}^3$ and 
$\nu_{\idx{S}} = 0.35$ for the Young's modulus, the density and the Poisson ratio, respectively.
The thin shell is discretized with tri-linearly interpolated hexahedral elements with EAS, comprising
two elements in thickness direction and $768$ elements over the surface. The internal pressure keeping 
the membrane inflated is set equal to the pre-shock state $p_{\idx{R}}$. Zero displacements in all three 
directions are prescribed for structural nodes located at the bottom of the shell at $y = 0\,\text{m}$. The time
integration factor $\theta = 0.5$ is chosen.

\begin{figure}
  \centering
  \includegraphics[width= 0.6\columnwidth]{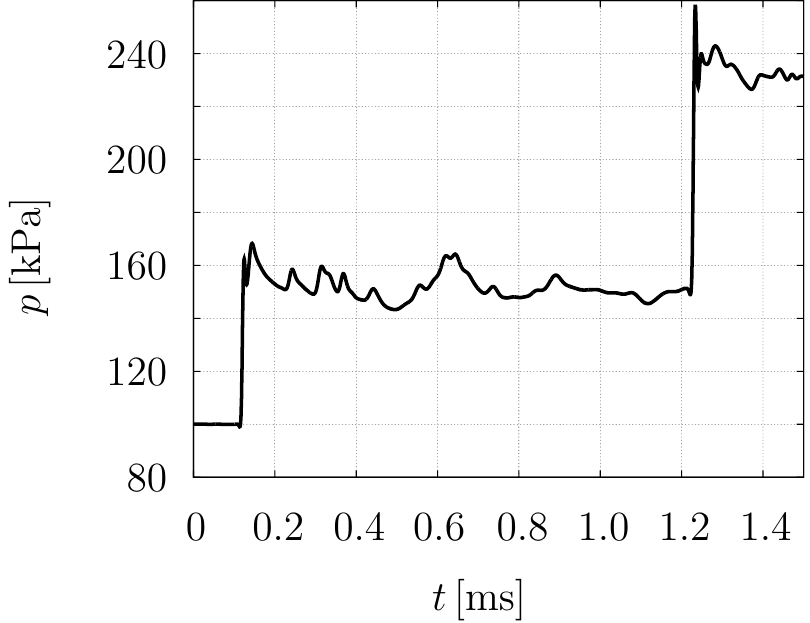}
  \caption{Pressure signal recorded at sensor position $(x,y,z) = (0,0.04,0)\,\text{m}$.}
  \label{fig:shock_membrane_pressure_sensor}
\end{figure}
\begin{figure}
  \centering
  \includegraphics[width= 1.0\columnwidth]{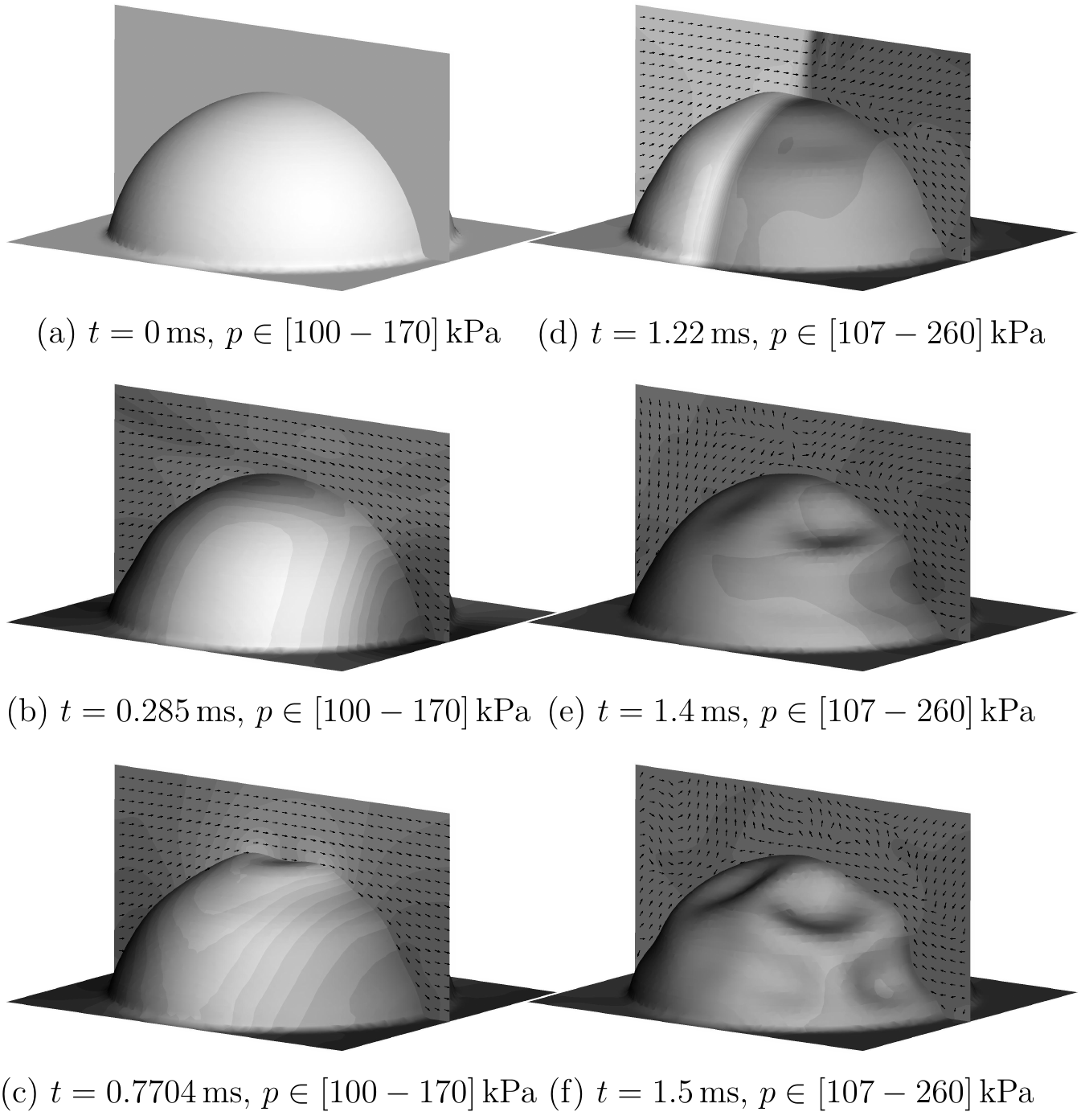}
  \caption{Pressure distribution together with uniform-length velocity vectors at different time 
instances. Every second vector is shown on the $xy$-midplane. Color scale ranges from 
white to black using $20$ equally spaced contour levels within the indicated pressure range.}
  \label{fig:shock_membrane_pressure_timeseries}  
\end{figure}
\begin{figure}
  \centering
  \includegraphics[width= 1.0\columnwidth]{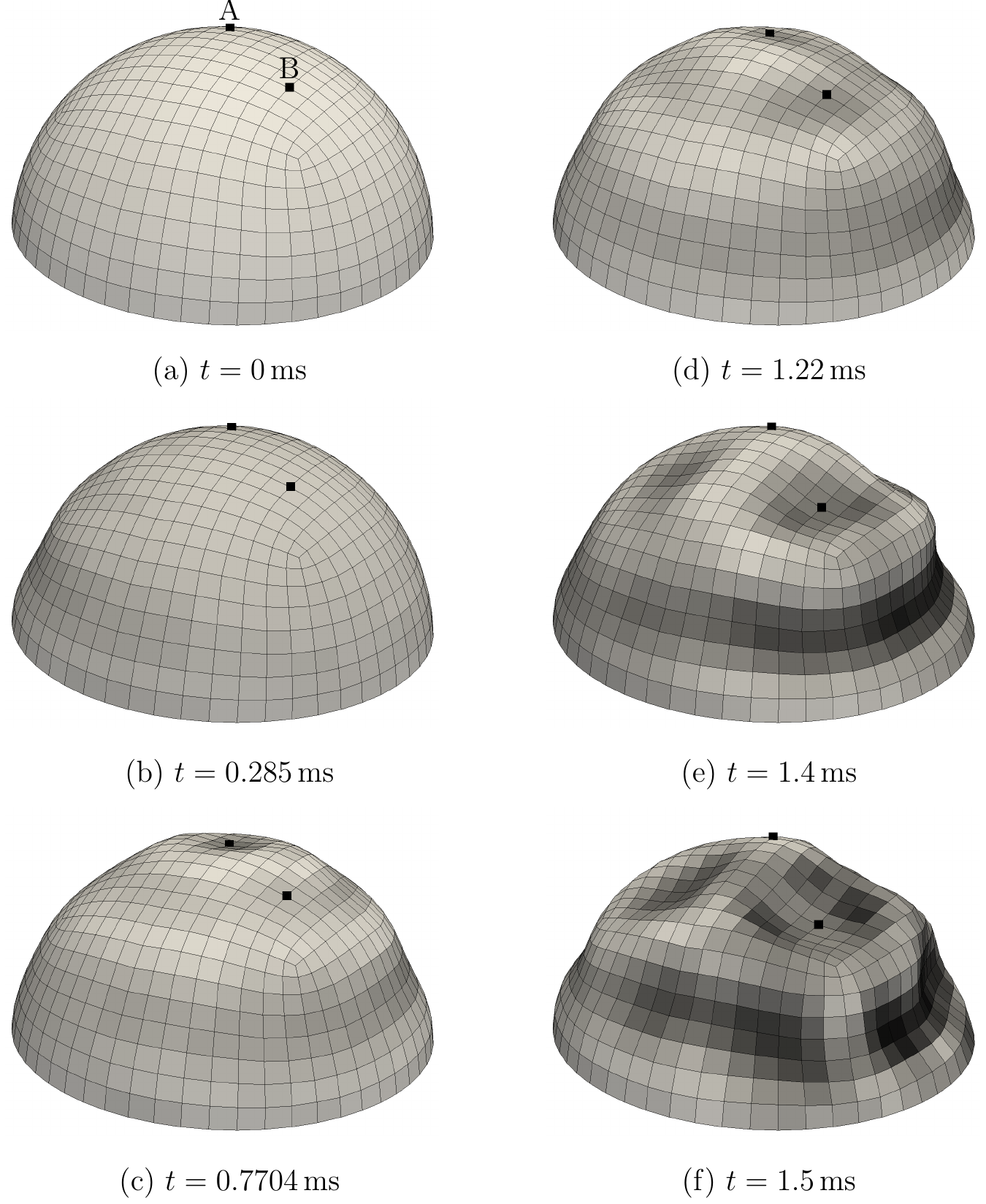}
  \caption{Norm of the Euler-Almansi strain tensor evaluated at each element center over
time. Color scale ranges from white to black using $26$ equally spaced contour levels for $\vert \tns{e} \vert_2
    \in [0-0.13]$. Monitoring points $A$ and $B$ are marked with squares.}
  \label{fig:shock_membrane_strain_timeseries}  
\end{figure}
\begin{figure}
  \centering
  \includegraphics[width= 0.7\columnwidth]{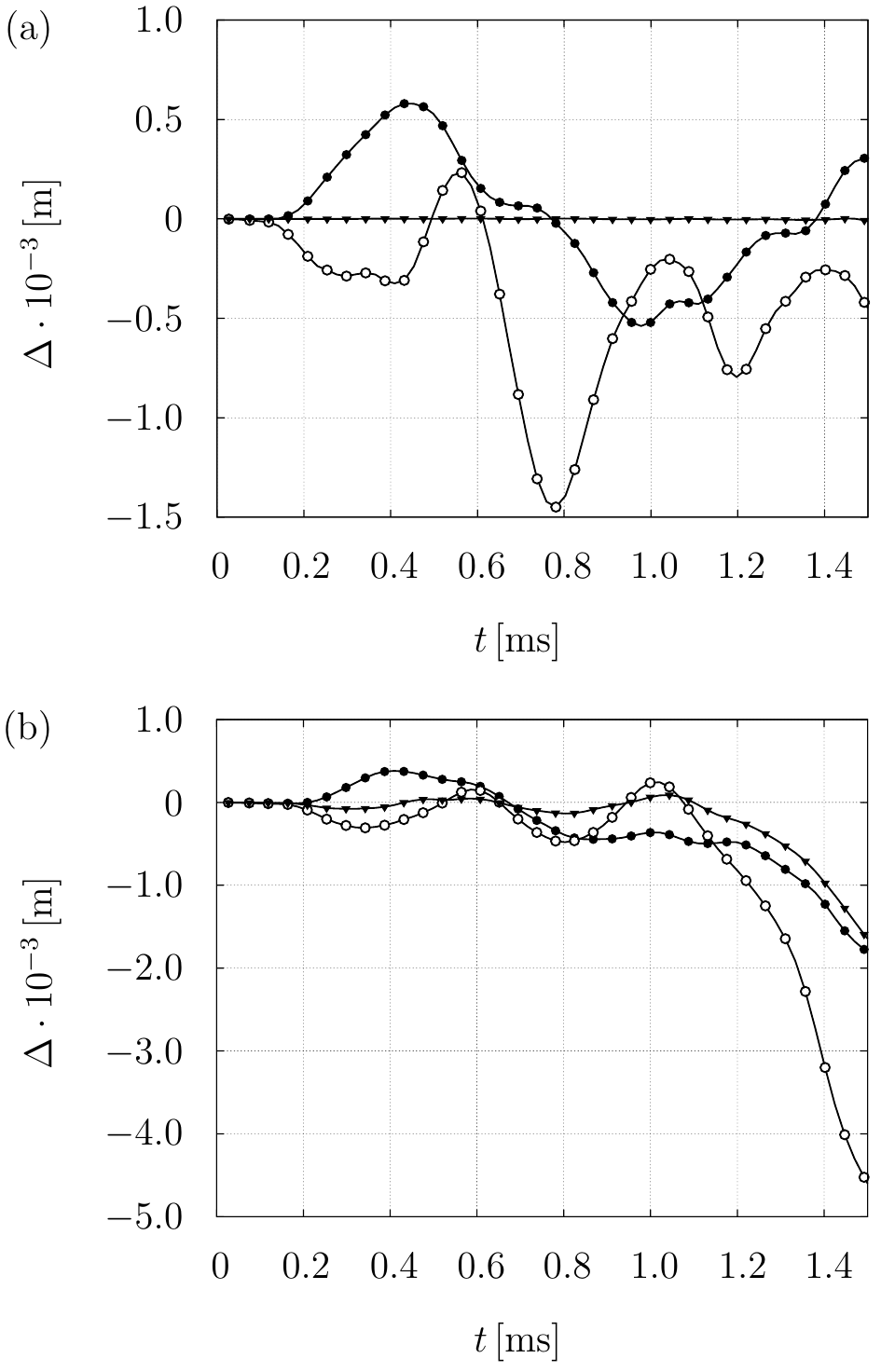}
  \caption{Time evolution of displacements at two different monitoring points (see
  Fig.~\ref{fig:shock_membrane_strain_timeseries}): (a) tip of membrane: monitoring point A, (b) monitoring point B. 
  (--- $\bullet$ ---) $\Delta x$, (--- $\circ$ ---) $\Delta y$, (--- $\blacktriangledown$ ---) $\Delta z$ }
  \label{fig:shock_membrane_displacements}
\end{figure}

Fig.~\ref{fig:shock_membrane_grid} depicts the computational domain and the fluid mesh in $xy$- and
$xz$-plane. In addition, we show the triangulated structural coupling interface, which is used for the cut
process in the fluid solver. Slip-wall boundary conditions are applied to all boundaries except for the inflow, where
all flow quantities are prescribed leading to a fully reflective boundary condition. In the region around the shell,
a uniform grid is used with cell sizes $\Delta x = \Delta y =
\Delta z = 0.001\,\text{m}$. In total, the fluid domain is discretized in space with $616,000$ cells. The time step
size is chosen to match a CFL number of $0.6$.

\begin{figure}
  \centering
  \includegraphics{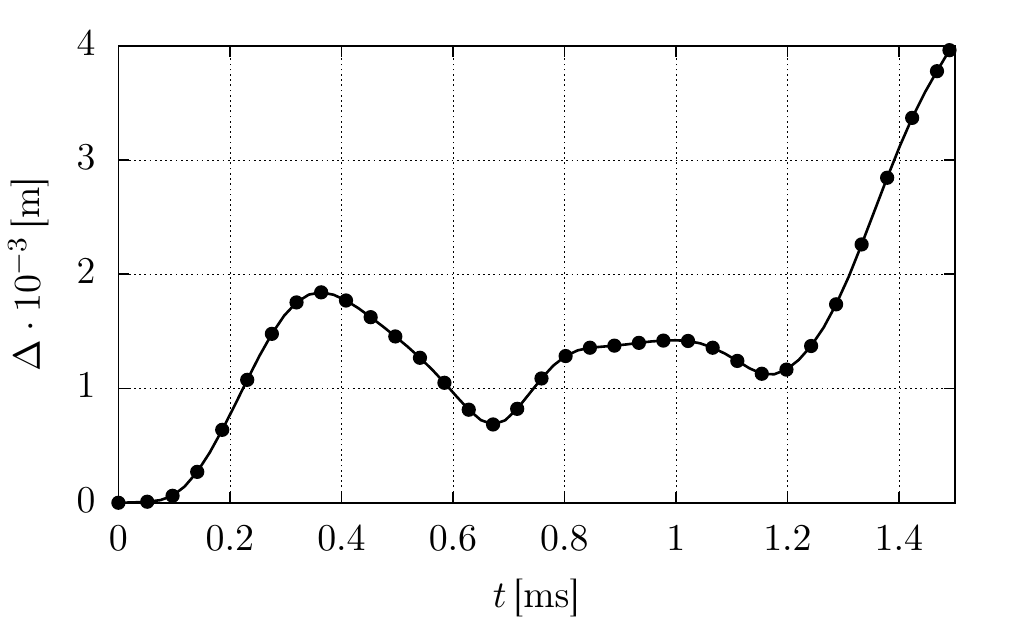}
  \caption{Time evolution of integral displacement magnitude.}
  \label{fig:shock_membrane_integral_disp}
\end{figure}

In Fig.~\ref{fig:shock_membrane_pressure_sensor}, the time evolution of the pressure signal recorded at the sensor 
position $P$ with $(x,y,z) = (0,0.04,0)\,\text{m}$ is shown. The jumps at approximately $t=0.1167\,\text{ms}$ and at
$t=1.2196\,\text{ms}$ mark the times when the shock wave passes the sensor.
Pressure distributions and velocity vectors at the $xy$-midplane are shown in
Fig.~\ref{fig:shock_membrane_pressure_timeseries} for different time instances. 
The corresponding strain distribution in the thin-walled shell is presented in
Fig.~\ref{fig:shock_membrane_strain_timeseries}. The norm of the
Euler-Almansi strain tensor $\vert \tns{e} \vert_2 = \sqrt{\tns{e} \, \colon \tns{e}}$ evaluated at each element center
of the top layer is chosen to illustrate the large deformations occurring during the buckling process.
Initially, the structure is undeformed and stress-free and the fluid is at 
rest, see Fig.~\ref{fig:shock_membrane_pressure_timeseries}(a) and Fig.~\ref{fig:shock_membrane_strain_timeseries}(a). 
Due to the
overpressure induced by the shock, Fig.~\ref{fig:shock_membrane_pressure_timeseries}(b), the windward side of the
membrane is compressed, see Fig.~\ref{fig:shock_membrane_strain_timeseries}(b), and is subsequently bouncing back due 
to its elastic behavior. At time $t = 0.7704\,\text{ms}$, buckling of the thin-walled shell occurs at its
tip, deflecting the flow as depicted in Fig.~\ref{fig:shock_membrane_pressure_timeseries}(c).
The displacement of the tip node at initial position $(x,y,z) = (0, 0.025, 0)\,\text{m}$ (monitoring point A) is given 
in
Fig.~\ref{fig:shock_membrane_displacements}(a): the $y$-deflection is approximately $1.5 \times 10^{-3}~\text{m}$ 
during this first
shock induced dimpling process.
As the shock hits the membrane after reflection at the end wall, 
see Fig.~\ref{fig:shock_membrane_pressure_timeseries}(d)-(f), the pressure increases again.
The membrane cannot sustain the additional load, and we observe the formation of buckling dimples, which are
symmetrically distributed with respect to the $xy$-midplane as shown
in Fig.~\ref{fig:shock_membrane_pressure_timeseries}(d)-(f) and Fig.~\ref{fig:shock_membrane_strain_timeseries}(d)-(f).
At $t = 1.5\,\text{ms}$, the norm of the Euler-Almansi strain in the most distorted regions rises up to
$0.127$, 
see Fig.~\ref{fig:shock_membrane_strain_timeseries}(f). Considering the monitoring point B, which is initially 
located at $(x,y,z) = (0.011912, 0.020912, 0.009308)~\text{m}$ in one of the dimples, a total 
deflection of $5.23 \times 10^{-3}~\text{m}$ is found, see Fig.~\ref{fig:shock_membrane_displacements}(b).

We refined the grids for both subdomains simultaneously and separately (not shown here for brevity) in 
order to reveal sensitivities with respect to the dynamic response of the thin-walled membrane. 
While the displacement of the membrane does not change significantly when varying the fluid resolution and keeping 
the structural discretization the same (maximum relative error of $2\%$ compared to a fluid grid with $\Delta x = 
\Delta y = \Delta z = 0.00025\,\text{m}$), we found that the dynamic response of the membrane and especially the 
occurring buckling mode can depend on the structural resolution. This observation confirms that 
buckling is highly sensitive with respect to imperfections of all kinds, including geometric imperfections 
\cite{Ramm2004}. Reliable prediction of buckling modes require realistic imperfection models, derived from the 
particular manufacturing process, to be included in the numerical model, which is beyond the scope of this paper. A 
well-defined quantity for such a configuration is the integral displacement magnitude shown in 
Fig.~\ref{fig:shock_membrane_integral_disp}. A maximum relative error of $3\%$ is found when comparing the 
present results to those of a four times finer mesh for both subdomains. A grid converged solution with respect to 
the integral displacement is obtained for a twice finer mesh.
Figure \ref{fig:shock_membrane_integral_disp} shows that the membrane starts to collapse at around $t = 1.2\,$ms, which 
coincides with the time when the shock wave, after reflection at the end wall, passes the pressure sensor $P$.

\section{Summary and conclusions}
\label{sec:conclusions}
The proposed finite volume -- finite element coupling approach for the interaction between a compressible fluid 
and a deformable structure is able to handle large and complex
three-dimensional deformations. We make use of a classical Dirichlet-Neumann partitioning in conjunction with a 
conventional serial staggered procedure for coupling of the two domains.

A representation of the interface within the fluid domain is achieved by means of a cut-element based IBM, which 
has been successfully extended to deformable structures for the first time. 
The presented framework leads to a non-matching discretization of the interface between both subdomains. A 
consistent data transfer has been established using a Mortar method, which preserves linear and angular momentum.
Piecewise constant ansatz 
functions are used for interpolating the fluid state as well as for the Lagrange multipliers 
on each single cut-element, allowing for a simple inversion of a diagonal matrix at negligible cost for the
evaluation of the discrete projection operator.
To the authors knowledge, this is the first time a cut-element method has been combined with a Mortar method for 
coupling the two subdomains in a consistent and efficient way.

The proposed coupling method has been validated through two-dimen\-si\-onal model problems involving rigid and 
deformable structures with large deformations. Our method correctly predicts the 
transient behavior of shock-loaded rigid and deformable structures. Moreover, good accuracy was achieved with 
respect to the correct prediction of flutter onset. The ability of our method to handle three-dimensional FSI problems 
involving large and complex structural deformations has been demonstrated through a newly proposed test case consisting 
of a flexible inflated thin shell interacting with a shock wave.

The current framework is limited to structures with a size larger than several fluid cells in order to fill the 
ghost-cell values properly. A remedy could be either an adaptive mesh refinement procedure for the flow solver 
or the decoupling of the ghost-cell method from the underlying Cartesian grid, which leads to additional degrees of 
freedom that need to be handled. In order to resolve the possibly different time-scales of both subdomains and increase 
the overall efficiency, subcycling should be considered for future work.

\section*{Acknowledgements}
The authors gratefully acknowledge support by the German Research Foundation (Deutsche Forschungsgemeinschaft) in the 
framework of the Collaborative Research Centre SFB/TRR 40 ``Fundamental Technologies for the Development of Future 
Space-Transport-System Components under High Thermal and Mechanical Loads''. Computational resources have been provided 
by the Leibniz Supercomputing Centre of the Bavarian Academy of Sciences and Humanities (LRZ).


\begin{appendix}
\section{Computational Performance}

\ctable[
        cap     = {Computational performance of the coupling framework for selected simulations.}, 
        caption = {Computational performance of the coupling framework for selected simulations.}, 
        label= {table:computational_costs},
        pos = {htbp},
        maxwidth = \columnwidth,
        captionskip = -0.4cm,
        doinside = \scriptsize,
        notespar,
        nosuper,
]{lrrrrrrr}{
	   \tnote[$a.$]{Only the rigid cylinder case is considered.}
           \tnote[$b.$]{Only the $50\,\text{mm}$ panel length case is considered.}
           \tnote[$c.$]{Only the $\text{Ma} = 2.3$ case is considered.}
        %
}{
\FL
Case (\# Run) &
$N^{ele}_{\idx{F}}$ &
$N^{ele}_{\idx{S}}$ &
$\frac{N^{ele}_{\idx{F}}}{N^{CPU}_{\idx{F}}}$& 
$\frac{N^{ele}_{\idx{S}}}{N^{CPU}_{\idx{S}}}$ & 
$T_{\idx{F}}\left( \%\right)$ &
$T_{\idx{S}}\left( \%\right)$ &
$T_{\idx{C}}\left( \%\right)$\\
\FL
\\
\textbf{\textsc{Cylinder}}\tmark[$a$]\\ \cmidrule(l){1-1}
$\#1$ & $2\cdot 10^3$ & $3.6\cdot 10^3$    & $2\cdot 10^3$ & $0.9\cdot 10^3$ & $1.7\,\%$ & 
$98.1\,\%$  & $0.2\,\%$ \\
$\#2$ & $8\cdot 10^3$ & $3.6\cdot 10^3$    & $8\cdot 10^3$ & $0.9\cdot 10^3$   & $4.4\,\%$  
& $95.4\,\%$ & $0.2\,\%$ \\
$\#3$ & $3.2\cdot 10^4$ & $3.6\cdot 10^3$  & $8\cdot 10^3$ & $0.9\cdot 10^3$   & $5.0\,\%$ 
& $94.8\,\%$ & $0.2\,\%$ \\
$\#4$ & $1.28\cdot 10^5$ & $3.6\cdot 10^3$ & $1.6\cdot 10^4$ & $0.9\cdot 10^3$ & $9.4\,\%$ 
& $90.4\,\%$ & $0.2\,\%$ \\
$\#5$ & $5.12\cdot 10^5$ & $3.6\cdot 10^3$ & $6.4\cdot 10^4$ & $0.9\cdot 10^3$ & $28.5\,\%$ 
& $71.3\,\%$ & $0.2\,\%$ \\
\\
\textbf{\textsc{Panel}}\tmark[$b$] \\ \cmidrule(l){1-1}
$\#1$ & $1.234\cdot 10^5$ & $1.3\cdot 10^2$    & $1.12\cdot 10^4$ & $4.3\cdot 10^1$ & $56.0\,\%$ 
&$43.4\,\%$  & $0.6\,\%$ \\
$\#2$ & $1.82\cdot 10^6$ & $1.3\cdot 10^2$    & $6.07\cdot 10^4$ & $6.5\cdot 10^1$ & $45.9\,\%$ & 
$53.9\,\%$  & $0.2\,\%$ \\
\\
\textbf{\textsc{Flutter}}\tmark[$c$] \\ \cmidrule(l){1-1}
$\#1$ & $1.65 \cdot 10^4$ & $1.6 \cdot 10^3$    & $4.125\cdot 10^3$ & $1.3\cdot 10^2$ & $13.9\,\%$ 
&$85.5\,\%$  & $0.6\,\%$ \\
\\
\textbf{\textsc{Membrane}} \\ \cmidrule(l){1-1}
$\#1$ & $6.16 \cdot 10^5$ & $1.536 \cdot 10^3$    & $3.08\cdot 10^4$ & $1.28\cdot 10^2$ & 
$45.2\,\%$ &$54.2\,\%$  & $0.6\,\%$
\LL}

The performance of the proposed coupling algorithm is summarized in Table~\ref{table:computational_costs}, 
where we show the percentage of time spent for the fluid solver $T_{\idx{F}}$, the structural solver $T_{\idx{S}}$,
and for the communication $T_{\idx{C}}$ between both codes for all considered test cases. $N^{ele}_{\idx{F,S}}$ 
represents the total number of elements used for the fluid and structural problem, respectively. 
$N^{ele}_{\idx{F,S}}/N^{CPU}_{\idx{F,S}}$ is the associated number of elements per CPU for each subdomain. The 
majority of the computational time is spent on advancing the solid domain, which, however, also includes load 
transfer with the Mortar method. Increasing the fluid resolution proportionally increases the number of 
cut-elements and thus the workload for the structural solver at the interface. Moreover, the implicit time
integration leads to an iterative solution procedure with at least two Newton iterations per coupling step to
obtain the solid state. The communication between both codes via Message Passing Interface typically
requires less than $1\,\%$ of the runtime. The current implementation of the staggered algorithm can be further 
optimized in terms of parallel efficiency. 
Furthermore, subcycling can significantly reduce computational cost of the structural solver and will be considered 
in future work.

\end{appendix}

\section*{References}
\bibliographystyle{elsarticle-harv}
\bibliography{fsi_paper}

\end{document}